\documentclass[a4paper,12pt]{article}

\usepackage{geometry}
\usepackage[english]{babel}
\usepackage[utf8]{inputenc}
\usepackage{tabularx}
\usepackage{amsmath}
\usepackage{amssymb}
\usepackage{amsthm}
\usepackage{enumitem}
\usepackage{mathtools}
\usepackage{bbold}
\usepackage{bm}
\usepackage{textcomp}
\usepackage{dsfont}
\usepackage{float}
\usepackage{algpseudocode}
\usepackage{color}
\usepackage{xcolor}
\usepackage{graphicx}
\usepackage{appendix}
\usepackage{subfig}
\usepackage{enumitem}
\usepackage{physics}
\usepackage{csquotes}
\usepackage[final]{hyperref}
\hypersetup{
	colorlinks=true,       
	linkcolor=black,       
	citecolor=blue,        
	filecolor=magenta,     
	urlcolor=blue
}

\usepackage[linesnumbered,ruled]{algorithm2e}
\makeatletter
\newcommand{\nosemic}{\renewcommand{\@endalgocfline}{\relax}}
\newcommand{\dosemic}{\renewcommand{\@endalgocfline}{\algocf@endline}}
\let\oldnl\nl
\newcommand{\nonl}{\renewcommand{\nl}{\let\nl\oldnl}}
\makeatother
\usepackage[backend=bibtex,doi=true,url=true,isbn=true]{biblatex}
\AtBeginBibliography{\footnotesize}
\theoremstyle{plain}
\newtheorem{theorem}{Theorem}[section]
\newtheorem{corollary}{Corollary}[theorem]
\newtheorem{lemma}[theorem]{Lemma}

\theoremstyle{definition}
\newtheorem{definition}{Definition}[section]
\newtheorem{remark}{Remark}[section]

\begin{document}

\title{Parameterizing Intersecting Surfaces via Invariants}

\author{Timon S. Gutleb\thanks{Department of Mathematics, University of British Columbia, Vancouver, V6T1Z2, BC, Canada}, Rhyan Barrett\thanks{Wilhelm-Ostwald-Institut für Physikalische und Theoretische Chemie, Universität Leipzig,  Germany}, 
Julia Westermayr\footnotemark[2] \thanks{ScaDS.AI (Center for Scalable Data Analytics and Artificial Intelligence) Dresden/Leipzig, Humboldtstraße 25, 04105 Leipzig}, 
Christoph Ortner\footnotemark[1]}
\date{\today}
\maketitle
\thispagestyle{empty}


\begin{abstract}
    We introduce and analyze numerical companion matrix methods for the reconstruction of hypersurfaces with crossings from smooth interpolants given unordered or, without loss of generality, value-sorted data. The problem is motivated by the desire to machine learn potential energy surfaces arising in molecular excited state computational chemistry applications. We present simplified models which reproduce the analytically predicted convergence and stability behaviors as well as two application-oriented numerical experiments: the electronic excited states of Graphene featuring Dirac conical cusps and energy surfaces corresponding to a sulfur dioxide ($SO_2$) molecule in different configurations. 
\end{abstract}

\section{Introduction}
Branching spaces, otherwise known as \emph{conical intersections}, are crucial for determining the dynamics of molecules and materials upon light excitation, molecular orbitals and material band structure. Their accurate description is especially important to advance our understanding of fundamental reactions like photosynthesis or vision \cite{mai_molecular_2020}. The relevant regions, however, are characterized by avoided crossings, i.e., two hypersurfaces that touch and become degenerate, but do not cross. This leads to non-smooth hyper-surfaces, making common machine learning models fail in their accurate modelling \cite{westermayr2020machine}.

The aim of this paper is thus to describe and analyze numerical methods for reconstructing multiple intersecting hyper-surfaces from unordered data using globally smooth interpolations. By \emph{unordered} we mean that the information ``which surface'' an evaluation point belongs to is unavailable. All surfaces throughout this work will be graphs. Our focus will be on theoretical aspects in a relatively idealized setting.

\begin{figure}\centering
     \subfloat[]
    {{ \centering \includegraphics[width=6.5cm]{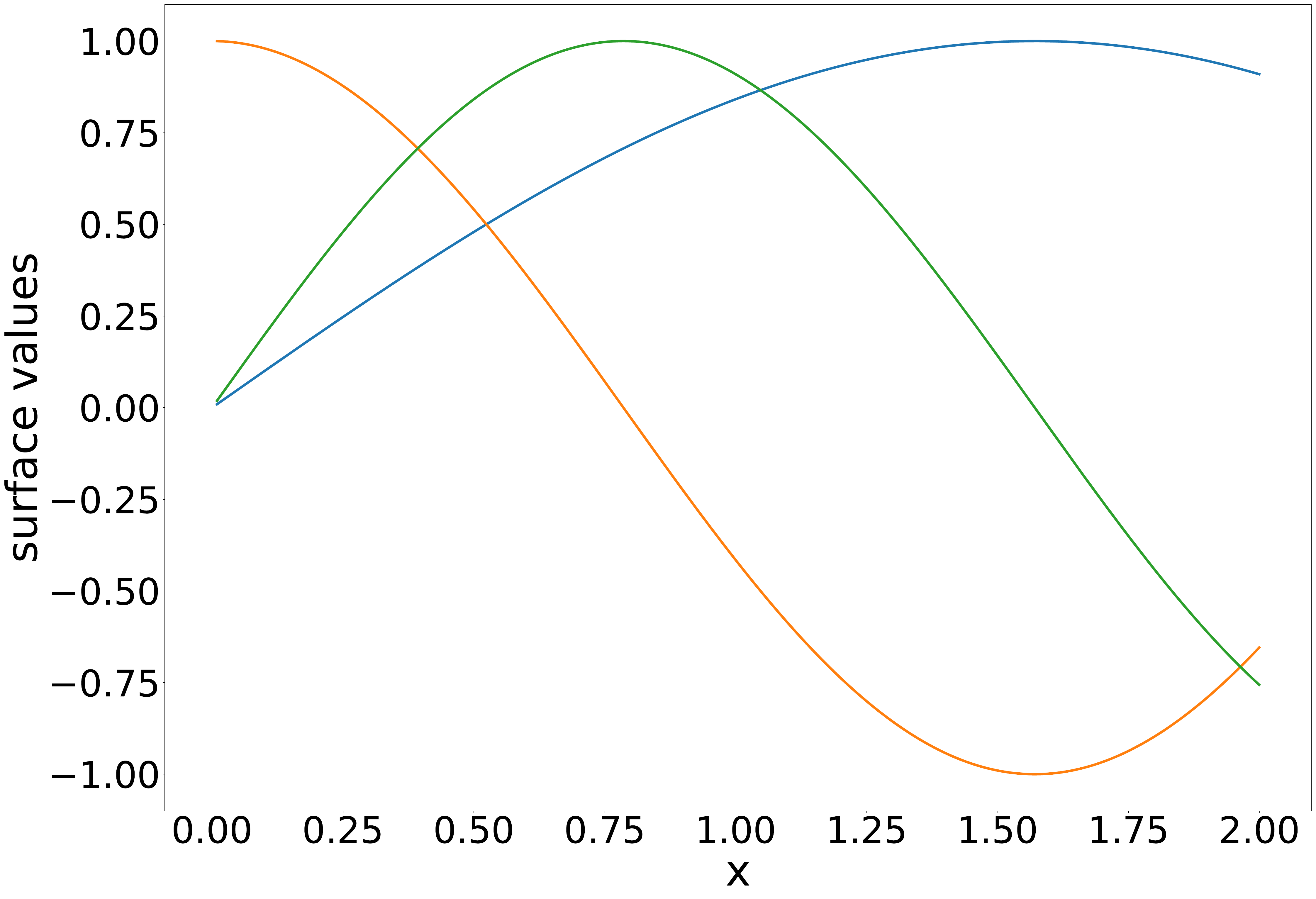} }}
     \subfloat[]
    {{ \centering \includegraphics[width=6.5cm]{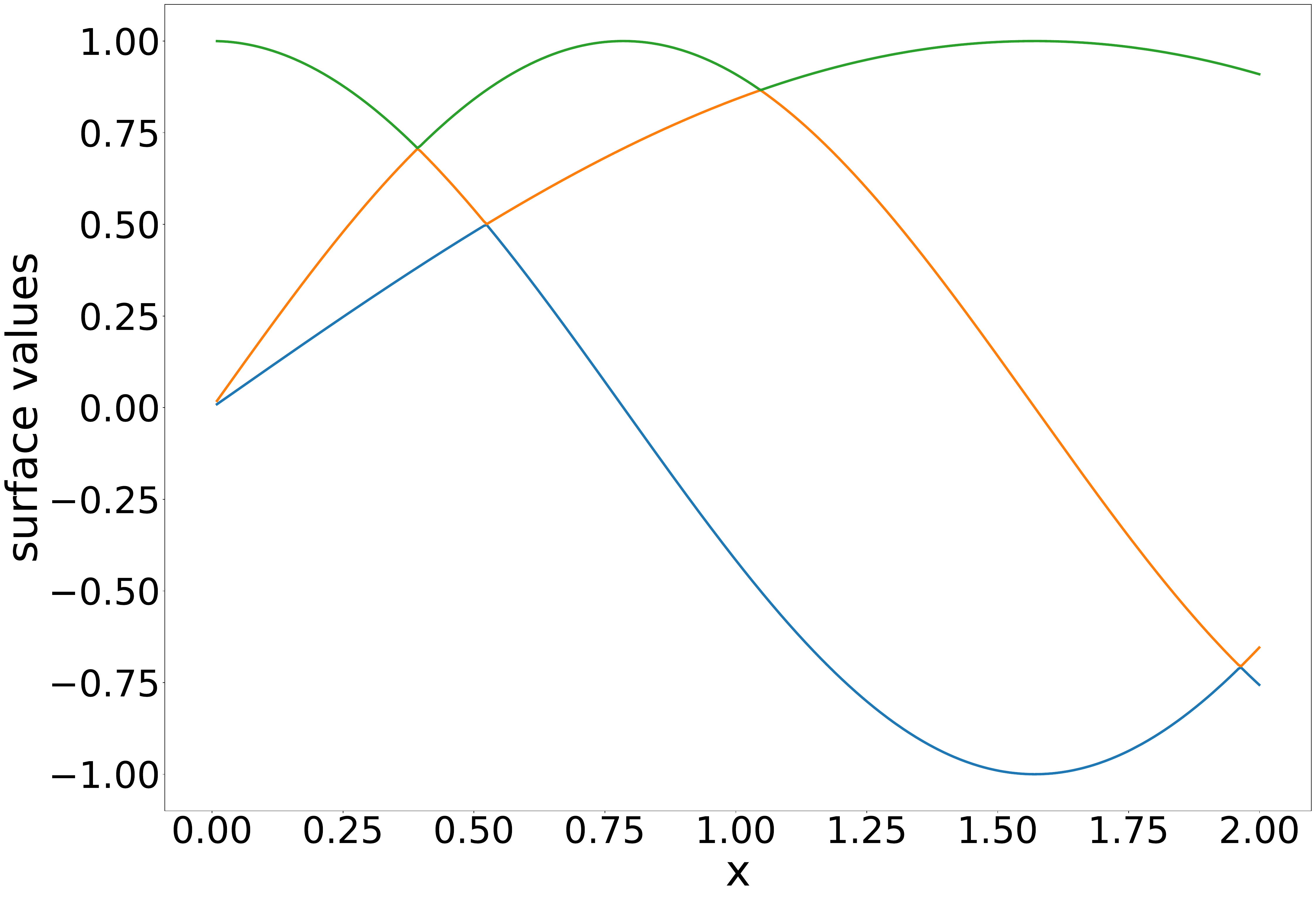} }}
    \caption{Intersecting, smooth sinusoidal curves in (a), sorted by value in (b), revealing non-smooth cusps in the ``$j$-th entry by value'' function.}
    \label{fig:firsttoyexample}
\end{figure}

By means of example, Figure \ref{fig:firsttoyexample}(a) shows three analytic curves, which are then value-sorted in (b) resulting in cusps at intersection points. The latter is, e.g., the form of the adiabatic potential energy surfaces produced in computational quantum chemistry. Well-known results on polynomial approximation, which we review in Sections \ref{sec:method} and \ref{sec:numericalanalysis}, tell us that the best-approximation rate in this case will be $O(1/n)$, where $n$ is the polynomial degree. On the other hand, if we knew the original analytic curves, i.e., if we could ``disentangle'' the intersections, then standard polynomial approximation schemes would guarantee an exponential or even super-exponential rate. 

Our investigation of reconstructing multiple intersecting graphs is motivated by applications in computational chemistry, where touching adiabatic potential energy surfaces of molecules are frequently studied using machine learning algorithms \cite{westermayr2020machine,westermayr_deep_2022,fedik_extending_2022,kulik_roadmap_2022,dral_molecular_2021,axelrod_excited_2022} relying on smooth interpolants but analogous challenges also arise for interpolating band diagrams in materials science applications \cite{sadat_machine_2020,xian_machine_2023,phani_dynamics_2017}. In most of these applications the surface values are given as a list of eigenvalues of an operator, with degenerate eigenvalues corresponding to the hypersurface crossings we are concerned with. The spectrum of this operator is typically expensive to compute, which motivates the need to construct computationally efficient surrogates. 

To circumvent the challenge posed by cusps near crossings of value-sorted surfaces, we introduce an approximation scheme based on approximating smooth invariants. For example, suppose the graphs are given by $f_i({\bf x}), i = 1, \dots, m,$ ${\bf x} \in \mathbb{R}^d$ then the power sum polynomials 
\[
    p_k({\bf x}) = \sum_{j = 1}^m f_{j}({\bf x})^k, \qquad k = 1, \dots, m
\]
are smooth functions of ${\bf x}$ (at least in the setting of Figure~\ref{fig:firsttoyexample}). Moreover, the power sum polynomials form a complete set of invariants:  from the values $(p_k)_{k=1}^m$ one can, in principle, reconstruct the surface values $(f_i)_{i = 1}^m$. By approximating the smooth $p_k$ we can recover a fast convergence rate under general and natural assumptions on the graphs $f_i$. Given an unseen input ${\bf x}$ one can then first predict the invariants $p_k({\bf x})$ from which approximations to the surfaces $f_i({\bf x})$ can then be reconstructed. 

In practice, this procedure is numerically subtle, and the purpose of this work is therefore to explore the numerical challenges that arise in its implementation. For example, we will work with a different set of invariants (the elementary symmetric polynomials) for which the reconstruction step can be formulated in terms of a companion matrix eigenvalue problem equivalent to root-finding for certain polynomials. In this setting there are several polynomial bases and companion matrix representations to choose from which are natural for different reasons and have distinct behavior in numerical practice. 

Machine learning seams of conical intersection for chemistry applications based on elementary symmetric polynomials was previously suggested by Opalka and Domcke as well as Schuurman, Neville and Wang in~\cite{opalka2013interpolation, wang2023machine} based on classical Frobenius companion matrices.
Our approach expands on these references in key aspects: First, we explore two natural alternative companion matrix approaches for this purpose -- Schmeisser companion and Chebyshev colleague matrices -- which have theoretical and practical advantages. In particular, using Frobenius companion matrices for root-finding problems is unstable in a lot of practically important cases. The second key difference is that we provide a numerical analysis description of these approaches. The latter is of particular importance for the algorithm's usefulness in the actual chemistry machine learning applications where understanding sensitivity of the methods to perturbations in the underlying data has critical consequences for performance and thus feasibility. We present numerous numerical experiments which agree with the results predicted by our analysis and test them in the application-relevant high noise regimes.

\section{Background}

\subsection{Elementary symmetric polynomials}
A {\em symmetric} function $f$ of $m$ variables is a function that is invariant under permutations of its arguments, i.e.,
\[
    f(x_1, \dots, x_m) = f(x_{\sigma 1}, \dots, x_{\sigma m}) \qquad \forall \sigma \in S_m,
\]
or $f \circ \sigma = f$ in short. Certain types of symmetric polynomials, such as of elementary, homogeneous, monomial and power sum types (cf. \cite[Ch.9]{loehr2017combinatorics}), are distinguished by various convenient properties and have received much attention in the fields of algebraic geometry, combinatorics and beyond \cite{egge2019introduction,Stanley_Fomin_1999,loehr2017combinatorics}. The numerical method we introduce in this paper relies primarily on properties of the elementary symmetric polynomials. 

\begin{definition}(Elementary symmetric polynomials, \cite[Ch.9]{loehr2017combinatorics} \& \cite[Ch.2]{egge2019introduction})\\
The $k$-th \emph{elementary symmetric polynomial} (ESP) in $m$ variables is defined by 
\begin{align*}
    s_k (x_1 , \ldots , x_m ):=\sum_{1\le  j_1 < j_2 < \cdots < j_k \le m} \hspace{2mm} \prod_{i=1}^k x_{j_i}.
\end{align*}
For example, 
\begin{align*}
s_0 (x_1, x_2, \dots,x_m) &= 1,\\
      s_1 (x_1, x_2, \dots,x_m) &= \sum_{1 \leq j \leq m} x_j,\\
  s_2 (x_1, x_2, \dots,x_m) &= \sum_{1 \leq j < k \leq m} x_j x_k,\\
  &\vdots\\
  s_m (x_1, x_2, \dots,x_m) &= \prod_{j=1}^m x_j
\end{align*}
and $s_k (x_1, x_2, \dots,x_m) = 0$ for $k>m$.
\end{definition}
Among the important properties of the elementary symmetric polynomials is that they establish a connection between the coefficients of a polynomial and its roots~\cite{vietae1646, girard1884, funkhouser1930short}, commonly known as Viète's formula. We refer to \cite[Section 3.2]{vinberg2003course} for a proof.  

\begin{theorem}[Viète's formula] \label{thm:viete} Let $p(x) = x^n + \sum_{j=0}^{n-1} a_j x^j$ be a degree $n$ monic polynomial with coefficients $a_j \in \mathbb{R}$ and (possibly complex and repeated) roots $r_j, j=1, \dots, n$, then 
\[
    s_k(r_1,\ldots,r_n) = (-1)^k a_{n-k}.
\]
\end{theorem}

\subsection{Companion matrices}\label{sec:Frobeniuscompanion}
In this section we review companion matrices of polynomials and some related constructions. Each square complex matrix $A$ has an associated characteristic polynomial given by $p_A(\lambda) = \det(\lambda I - A)$, whose roots are precisely the eigenvalues of $A$ \cite{strang2022introduction,horn2012matrix}. Companion matrices~\cite{Frobenius1879, loewy_begleitmatrizen_1920, mac_duffee_theory_1933} arise as the answer to the converse question: Given a polynomial $p(x)$, can one construct a matrix $A_p$ whose eigenvalues are the roots, i.e. $p(x) = \det(x I - A_p)$? This question is particularly interesting for numerical computations since it allows one to replace root finding problems with well-understood eigenvalue problems  \cite{trefethen2011six,aurentz2015fast,aurentz2018fast,edelman1995polynomial}.

\begin{theorem}(Frobenius companion matrix, \cite[Ch.3]{horn2012matrix}) \label{def:frobeniuscompanion} Let $p(x) = x^n + \sum_{k=0}^{n-1} a_{k} x^{k}$ be a degree $n$ monic polynomial. We define the \emph{companion matrix} $A_p^{\rm F}$ of $p(x)$ as \begin{align*}
   A_p^{\rm F} =\begin{pmatrix}
0 & 0 & \dots & 0 & -a_0 \\
1 & 0 & \dots & 0 & -a_1 \\
0 & 1 & \dots & 0 & -a_2 \\
\vdots & \vdots & \ddots & \vdots & \vdots \\
0 & 0 & \dots & 1 & -a_{n-1}
\end{pmatrix}.
\end{align*}
Then, the eigenvalues of the companion matrix $A^{\mathrm{F}}_p$ are the roots of $p(x)$.
\end{theorem}

Since matrix eigenvalues are preserved under similarity transforms, the Frobenius companion matrix $A^{\mathrm{F}}_p$ is not a unique solution to $p(x) = \det(x I - A)$; see, e.g., the pentadiagonal Fiedler companion matrices \cite{fiedler2003note, mackey2013continuing, fiedler1990expressing}. Our work 
uses Schmeisser's related construction of symmetric tridiagonal companian matrices \cite{schmeisser1993real}. 

\begin{theorem}[Schmeisser companion matrix, \cite{schmeisser1993real}]\label{thm:schmeisser} 
    Let $p(x) = x^n + \sum_{k=0}^{n-1} a_{k} x^{k}$,  $a_j \in \mathbb{R}$ be a monic polynomial of degree $n$ with only real (but not necessarily distinct) roots. Then there exists a real symmetric tridiagonal matrix $A^{\mathrm{S}}_p$ with non-negative off-diagonal entries
\begin{align}
    A^{\rm S}_p = \begin{pmatrix}
        -q_1(0) && \sqrt{c_1} && 0 && \hdots && 0\\
        \sqrt{c_1} && -q_2(0) && \sqrt{c_2} && \ddots && 0\\
        0 && \sqrt{c_2} && -q_3(0) && \ddots && \vdots \\
        \vdots && \ddots && \ddots && \ddots && \sqrt{c_{n-1}}\\
        0 && \hdots && 0 && \sqrt{c_{n-1}} && -q_n(0)
    \end{pmatrix},
\end{align}
    for which $p(x)$ is the characteristic polynomial, that is
    \begin{align*}
        p(x) = \det(xI-A_p^{\mathrm{S}}).
    \end{align*}
    Furthermore, Algorithm~\ref{alg:constructiveschmeisser} provides a constructive method to obtain the entries of $A^{\mathrm{S}}_p$.
\end{theorem}

\begin{algorithm}
    \caption{Construction of sym. tridiagonal Schmeisser companion matrix \cite{schmeisser1993real}}\label{alg:constructiveschmeisser}\vspace{2mm}
    \KwIn{Monic polynomial $p(x) = x^n + a_{n-1}x^{n-1} + ... a_0$,  $a_j \in \mathbb{R}$, with exclusively real roots of arbitrary multiplicity. Coefficient vectors sorted by ascending order, i.e. $\mathbf{a} = (a_0,...a_{n-1},1)$ and indexed by $1:\text{end}$.}
    \KwOut{Diagonal and off-diagonal bands $(d,s)$ of $A^{\mathrm{S}}_p$ in Theorem \ref{thm:schmeisser}.} 
        Set $y_1 \gets \mathbf{a}$ and $y_2
        \gets \frac{\mathbf{a}'}{n}$ where $\mathbf{a}'$ are the coefficients of $p'(x)$.\\
        \For{$k=1$ \KwTo $n$}{
     $(q,r) \gets$ Quotient and remainder of polynomial division of $y_1$ by $y_2$.\\
 \If{$k < n$}
{
    $y_1 \gets y_2$\\
    $y_2 \gets \frac{r}{r[\text{end}]}$\\
    $s[k] \gets -r[\text{end}]$
}
$d[k] \gets -q[1]$
    }
        
    \Return $(d,s)$
\end{algorithm}

A key advantage of Schmeisser's construction is that the eigenvalues of symmetric tridiagonal matrices are not only guaranteed to be real but there are also efficient and well-understood numerical algorithms for the computation of these eigenvalues which guarantee real results, cf. \cite{parlett1998symmetric,golub1962bounds}. 

\subsection{Colleague matrices and Chebyshev polynomials}\label{sec:chebintro}
The final concept to review are the \emph{colleague} matrices, which first require a brief discussion of root-finding with Chebyshev polynomials.
It is a classical question how polynomials behave under perturbations of their coefficients \cite{gautschi1972condition,gautschi1979condition,wilkinson1959evaluation1,wilkinson1959evaluation2,mosier1986root}. A first cautionary tale is immediately found in $x^2$, where perturbation to $x^2+\epsilon$ by a small $0 < \epsilon \ll 1$ changes the roots from a multiplicity two real root at $0$ to entirely imaginary roots $\pm i \sqrt{\epsilon}$ with the error scaling with $\sqrt{\epsilon} \gg \epsilon$. In this case the multiplicity causes the issue but as was shown by Wilkinson in \cite{wilkinson2023rounding} this is not always the case: Wilkinson's example polynomial $\prod_{n=1}^{20} (x-n)$ has exclusively real and distinct roots $\{1,2,...,20\}$ but when expanded into a monomial series has extremely poor numerical root-finding conditioning, cf. \cite{mosier1986root} and \cite[Part 3]{trefethen2022numerical}.

Despite these observations, polynomial root-finding is not necessarily ill-conditioned. As described in \cite{trefethen2011six}, using a companion matrix based root-finding method for a polynomial expressed in the monomial basis is well-behaved if the roots are near the complex unit circle but otherwise becomes ill-conditioned~\cite{trefethen2019approximation,sitton2003factoring}. In general, the conditioning with respect to perturbed coefficients depends on (1) the multiplicities of the roots, (2) the basis in which a polynomial is expressed and (3) the location of its roots in the complex plane. 

The Chebyshev polynomials of first kind $\{T_j(x)\}_{j\in\mathbb{N}_0}$ are a basis of polynomials orthogonal with respect to the following inner product:
\begin{align*}
    \langle T_n, T_m \rangle_T = \int_{-1}^1 T_n(x)\,T_m(x)\,\frac{1}{\sqrt{1-x^2}} \mathrm{d}x = 
\begin{cases}
0             & ~\text{ if }~ n \neq m, \\
\pi           & ~\text{ if }~ n=m=0, \\
\frac{\pi}{2} & ~\text{ if }~ n=m \neq 0
\end{cases}
\end{align*}
and have many appealing theoretical and numerical properties \cite{rivlin_chebyshev_2020,trefethen2019approximation,DLMF}. Importantly, the roots of $p(x) = \sum_{j=0}^n b_j T_j(x)$ are well-conditioned as a function of $b_j$ if its roots lie on or near $[-1,1]$ and an eigensolver for the matrix introduced in the following theorem is used \cite{specht1956lage,good1961colleague,trefethen2019approximation,serkh_provably_2021,noferini2017chebyshev,nakatsukasa2016stability}. 

\begin{theorem}[Colleague matrix, \cite{good1961colleague,trefethen2019approximation}] \label{thm:colleaguedef}
    The roots of the polynomial $p(x) = T_n(x) + \sum_{j=0}^{n-1} b_j T_j(x)$,  $b_j \in \mathbb{R}$, are the eigenvalues of the following matrix which is commonly called the \emph{colleague matrix}:
  \begin{align}
A^{\mathrm{C}}_p = H - \frac{1}{2}e_1 c^\top = \left(
  \begin{array}{ccccc}
0 & \frac{1}{2} &  &  &   \\
\frac{1}{2} & 0 & \frac{1}{2}  &   \\
& \frac{1}{2} & \ddots & \ddots & \\
&  & \ddots & 0 & \frac{\sqrt{2}}{2} \\
&  &  & \frac{\sqrt{2}}{2} & 0
  \end{array} \right)
- \frac{1}{2} e_1 \left( 
  \begin{array}{ccccc}
b_{n-1} & \cdots & b_1 & \sqrt{2} b_0   
  \end{array} \right),
    \label{colleague}
  \end{align}
\end{theorem}

\begin{remark} The colleague matrix $A^{\mathrm{C}}_p$ in Theorem \ref{thm:colleaguedef} is an upper Hessenberg matrix and is expressed as a sum of a real symmetric matrix and a rank-1 perturbation.

An overview of numerical algorithms, as well as a provably component-wise backward stable $O(n^2)$ QR algorithm for the diagonalization of colleague matrices was recently given by Serkh and Rokhlin \cite{serkh_provably_2021}. 
\end{remark}

\section{Reconstruction of Multi-Surfaces} \label{sec:method}
\subsection{Problem statement and preliminaries}
We assume that a union of $m \geq 2$ (hyper-)surfaces is given as graphs over some domain $\Omega \subseteq \mathbb{R}^{d}$ and call such a collection a \emph{multi-surface} for simplicity. Corresponding with our applications, we assume that at some finite collection of points $x\in\Omega \subseteq \mathbb{R}^{d}$ we can evaluate the surfaces sorted by their values. Any ordering can always be sorted to match this scenario and thus this comes at no loss of generality. Notably this means that there will generally be non-smooth cusps in the $j$-th entry function of our data as shown in Figure \ref{fig:firsttoyexample} causing challenges for globally smooth, e.g. polynomial, approximation approaches.

Our aim is to reconstruct the multi-surface from given value-ordered point-wise data using a globally smooth interpolation method with good numerical properties (convergence, error control, etc.). This goal is motivated in large part by the use of machine learning approaches in the intended applications which rely on learning globally smooth interpolations. The following theorem provides the rigorous justification for reconstructing the multi-surfaces in this way as opposed to simply approximating the non-smooth surfaces featuring cusps, cf. Figure \ref{fig:firsttoyexample}, directly.
\begin{theorem}\label{thm:chebconvergence}\cite[Theorem 7.2]{trefethen2019approximation}
For $\nu \in \mathbb{N}$, let a function $f$ and its derivatives up to $f^{(\nu-1)}$ be absolutely continuous on $[-1,1]$ and let $f^{(\nu)}$ satisfy $V = \|f^{(\nu+1)}\|_1 < \infty$ (bounded variation). Then for any $n > \nu$, the degree $n$ Chebyshev interpolation of $f$ denoted $p_n$ satisfies
\begin{align*}
\|f - p_n \|_\infty \leq \frac{4V}{\nu (n-\nu)^\nu} = O(n^{-\nu}).
\end{align*}
If $f$ is analytic in $[-1, 1]$ then there exist constants $c, \alpha > 0$ such that 
\[
    \| f - p_n \|_\infty \leq c e^{- \alpha n}. 
\]
\end{theorem}
The cusps seen in Figure \ref{fig:firsttoyexample} correspond to the case $\nu = 1$ in Theorem \ref{thm:chebconvergence} (when rescaled to $[-1,1]$), i.e. the cusp-containing surfaces are absolutely continuous but their derivatives which exist almost everywhere are not (though they have bounded variation). As a result, one should not expect convergence better than $O(n^{-1})$ when using a global Chebyshev approximation for the $j$-th entry multi-surface functions. One readily observes that if the ESPs are smooth then a Chebyshev interpolation of these functions converges exponentially in degree (interpolation point spacing and number, among other considerations, may put practical limitations on this theoretical behavior).

In the next sections we describe three closely related methods for the reconstruction of multi-surfaces from global interpolants of invariants. The three methods correspond, respectively, to employing Frobenius companion matrices, Schmeisser companion matrices and Chebyshev colleague matrices in the reconstruction, each with distinct conceptual or numerical advantages. 

\subsection{Reconstruction from invariants, Frobenius variant}
We begin with the conceptually simplest variant using the classical Frobenius companion matrices and motivate the other variants as modifications thereof.
We first describe the idea of the method for the simple case of three intersecting, value-sorted surfaces $\mathbf{f}(\mathbf{x}) = \text{sort}(f_1(\mathbf{x}),f_2(\mathbf{x}), f_3(\mathbf{x}))$, then present the general case methods. First, we construct the ESP values of the surfaces at each given point, i.e. for the case of three surfaces we have
\begin{align*}
    s_1(\mathbf{f}(\mathbf{x})) &=  \mathbf{f}(\mathbf{x})_1+\mathbf{f}(\mathbf{x})_2+\mathbf{f}(\mathbf{x})_3,\\
    s_2(\mathbf{f}(\mathbf{x})) &= \mathbf{f}(\mathbf{x})_1 \mathbf{f}(\mathbf{x})_2 + \mathbf{f}(\mathbf{x})_1 \mathbf{f}(\mathbf{x})_3 + \mathbf{f}(\mathbf{x})_2 \mathbf{f}(\mathbf{x})_3 ,\\
    s_3(\mathbf{f}(\mathbf{x})) &= \mathbf{f}(\mathbf{x})_1 \mathbf{f}(\mathbf{x})_2 \mathbf{f}(\mathbf{x})_3.
\end{align*}
In a typical application scenario where one cannot evaluate the surfaces at arbitrary points but must instead work with pre-existing data, one would fit a (typically global) approximation of the three ESPs using e.g. orthogonal polynomials. 
Note that the $j$-th element of $\mathbf{f}(\mathbf{x})$, denoted $\mathbf{f}(\mathbf{x})_j$, is generally \emph{not} equal to $f_j(\mathbf{x})$ due to the lack of sorting, cf. Figure \ref{fig:firsttoyexample}. However, by construction the ESPs are invariant under permutation of their arguments, meaning that any information about the sorting applied to $\mathbf{f}(\mathbf{x})$ is lost in this process. By Viète' formula (Theorem \ref{thm:viete}), the values of $\mathbf{f}(\mathbf{x})_j$ can be recovered from the ESPs by finding the roots of the monic polynomial
\begin{align*}
    p(y) = y^3 - s_1(\mathbf{f}(\mathbf{x})) y^2 + s_2(\mathbf{f}(\mathbf{x})) y - s_3(\mathbf{f}(\mathbf{x})),
\end{align*}
which we can achieve by computing the eigenvalues of its $3\times 3$ Frobenius companion matrix $A^{\mathrm{F}}_p(\mathbf{x})$ as defined in Section \ref{sec:Frobeniuscompanion} at each point $\mathbf{x} \in \Omega$, with
\begin{align*}
    A^{\mathrm{F}}_p(\mathbf{x}) = \begin{pmatrix}
        0 & 0  & s_3(\mathbf{f}(\mathbf{x}))\\
        1 & 0  & -s_2(\mathbf{f}(\mathbf{x}))\\
        0 & 1  & s_1(\mathbf{f}(\mathbf{x}))\\
    \end{pmatrix}.
\end{align*}
The order of the returned eigenvalues may differ depending on the solver, so to guarantee consistency we sort the obtained eigenvalues $\lambda_{A^{\mathrm{F}}_p(\mathbf{x})}$ by value to obtain
\begin{align*}
    \text{sort}(\lambda_{A_p(\mathbf{x})}) = \mathbf{f}(\mathbf{x}) = \text{sort}(f_1(\mathbf{x}),f_2(\mathbf{x}), f_3(\mathbf{x})).
\end{align*}

In describing this method for the reconstruction of $\mathbf{f}(\mathbf{x})$ from its ESPs, it may appear as if we have walked in a circle and simply regained what we started from: We were given the value-sorted $\mathbf{f}(\mathbf{x})$, constructed the ESPs from it and recovered the value-sorted $\mathbf{f}(\mathbf{x})$ from the associated companion matrix eigenvalue problem. The important observation which turns this into a useful framework is that while $\mathbf{f}(\mathbf{x})$ has cusps and is thus not well-approximated by globally smooth approximations (e.g. orthogonal polynomials), the ESPs are guaranteed to be smooth if \emph{any} underlying smooth ordering of $\mathbf{f}(\mathbf{x})$ exists (by permutation invariance and linearity). Furthermore, as we will explore for conical cusps in the numerical experiments in Section~\ref{sec:numericaltoy2cusp}, the only condition for well-approximable ESPs is their own smoothness which may be given even in situations where the underlying surfaces are not smooth. 

The general procedure for arbitrary number of surfaces is described in Algorithm~\ref{alg:frobeniusmethod}, which is exact in exact arithmetic. As mentioned in the introduction, Algorithm~\ref{alg:frobeniusmethod} is conceptually closely related to ideas suggested in \cite{wang2023machine, opalka2013interpolation}.

\begin{algorithm}
    \caption{Pointwise reconstruction from smooth interpolants via Frobenius companion matrix}\label{alg:frobeniusmethod}\vspace{2mm}
    \KwIn{Elementary symmetric polynomials $\mathbf{s}(\mathbf{f}(\mathbf{x}))$ corresponding to an underlying multi-surface $\mathbf{f}(\mathbf{x})$ with $m$ entries at a given point $\mathbf{x}' \in \Omega \subset \mathbb{R}^d$.}
    \KwOut{Pointwise values of multi-surface $\mathbf{f}(\mathbf{x}')$ sorted by value.} 
       $A^{\mathrm{F}}_p \gets$ Companion matrix with coefficients $(-1)^{m+1}s_m(\mathbf{f}(\mathbf{x}))$ via Thm.~\ref{def:frobeniuscompanion}.\\
       $ \lambda \gets \text{eigenvalues}(A^{\mathrm{F}}_p)$
        
    \Return $\text{sort}(\lambda)$
\end{algorithm}

There are several drawbacks of the method presented in Algorithm \ref{alg:frobeniusmethod}, first and foremost the fact that near cusps (i.e. multiple roots of the underlying polynomial) the eigenvalue computation can produce incorrect complex roots due to numerical errors, which when interpreted in terms of their real or absolute values cause surfaces to clamp together instead of intersecting. Next, we discuss the Schmeisser variant which helps address some of the shortcomings.

\subsection{Symmetric tridiagonal Schmeisser variant}
In Algorithm \ref{alg:schmeissermethod} we present an algorithm which reconstructs pointwise multi-surface values from the globally smooth ESPs thereof using the Schmeisser companion matrix which due to its real symmetric nature is guaranteed to have real eigenvalues. If appropriate eigensolvers are used (denoted in Algorithm \ref{alg:schmeissermethod} as \emph{eigh()} due to the NumPy convention), the numerical eigensolver is also guaranteed to produce real eigenvalues, making complex roots impossible -- albeit with a caveat.

\begin{algorithm}
    \caption{Pointwise reconstruction from smooth interpolants from the Schmeisser companion matrix}\label{alg:schmeissermethod}\vspace{2mm}
    \KwIn{Elementary symmetric polynomials $\mathbf{s}(\mathbf{f}(\mathbf{x}))$ corresponding to an underlying multi-surface $\mathbf{f}(\mathbf{x})$ with $m$ entries at a given point $\mathbf{x}' \in \Omega \subset \mathbb{R}^d$.}
    \KwOut{Pointwise values of multi-surface $\mathbf{f}(\mathbf{x}')$ sorted by value (guaranteed real).} 
       $A^{\mathrm{S}}_p \gets$ Schmeisser companion matrix for $(-1)^{m+1}s_m(\mathbf{f}(\mathbf{x}))$ via Algorithm \ref{alg:constructiveschmeisser}\\
       $ \lambda \gets \text{eigh}(A^{\mathrm{S}}_p)$
        
    \Return $\text{sort}(\lambda)$
\end{algorithm}

Analogous to the discussion in the previous section, Algorithm \ref{alg:schmeissermethod} is exact in exact arithmetic. However, if there is sufficient noise on the input, then it is possible that the polynomial with theoretically guaranteed real roots is perturbed into a polynomial with complex roots (cf. the discussion in Section \ref{sec:chebintro} for monomial series). If the Schmeisser construction in Algorithm \ref{alg:schmeissermethod} is attempted with a polynomial whose roots are complex, then it will fail since a polynomial with complex roots cannot have a Hermitian companion matrix. In practice the failing occurs such that the $c_j$ on the off-diagonals of the Schmeisser matrix (see the definition in Theorem \ref{thm:schmeisser}) become negative due to noise or numerical round-off, causing the square root to be complex, resulting in a non-Hermitian matrix if complex square roots are allowed. 

One way to avoid this is motivated by a classical result: Again analogous to the Frobenius method, the surface crossings occur precisely where two or more eigenvalues are identical -- more precisely, no surfaces cross if all eigenvalues of the companion matrix are \emph{simple}. Unlike for Frobenius companion matrices, however, one can immediately determine based on the off-diagonal entries of the Schmeisser companion matrix whether two or more surfaces are close to each other or crossing due the following classical result for symmetric tridiagonal matrix eigenvalues:

\begin{lemma}(\cite[Section 7]{parlett1998symmetric}) Let $A$ be a real symmetric tridiagonal matrix. If all the off-diagonal entries are positive, then all eigenvalues of $A$ are simple.
\end{lemma}
Note that this implies that crossings occur exactly where an off-diagonal term $c_j$ of the Schmeisser companion matrix vanishes. To avoid small negative values in the $c_j$ from erroring, one can introduce a filtering step into Algorithm \ref{alg:schmeissermethod}, e.g. $c_j \gets \mathrm{max}(0,c_j)$, effectively corresponding to a manual clamping of values. In contrast to the Frobenius method, however, (1) proximity to cusps can be easily read off of the smallest values on the off-diagonals and thus (2) the filtering step can be done in many different, controlled ways and (3) the final results if filtered appropriately are always guaranteed to be real.

\subsection{Chebyshev Colleague matrix variant}
A Chebyshev colleague matrix variant of Algorithms \ref{alg:frobeniusmethod} and \ref{alg:schmeissermethod} is straightforward to construct as seen in Algorithm \ref{alg:chebmethod} but requires access to Chebyshev coefficients $b_j$ such that $$T_m(y)+\sum_{j=0}^{m-1} b_j(\mathbf{x}') T_j(y) = \sum_{k=0}^{m} (-1)^{m-k} s_{m-k}(\mathbf{f}(\mathbf{x}')) y^k .$$
Generically, the Chebyshev coefficients of any sufficiently regular function $f$ can be obtained by the Chebyshev inner product (cf. \cite[Ch.4]{lanczos_applied_1988})
\begin{align*}
    b_j = \frac{2}{\pi} \langle f, T_j \rangle_T =  \int_{-1}^1 f(x) \,T_j(x)\,\frac{1}{\sqrt{1-x^2}} \mathrm{d}x = \frac{2}{\pi} \int_0^\pi f(\cos(\theta)) \cos(j \theta) \mathrm{d}\theta,
\end{align*}
where such an integral can be computed numerically or analytically.
To our knowledge there is no closed-form representation of the required ESP-as-Chebyshev coefficients but this poses no difficulty for our proposed algorithm as for any fixed number of intersecting surfaces $m$ the ESP-as-Chebyshev coefficients can be symbolically pre-computed using the following monomial conversion rule often used in the context of power sum economization \cite{lanczos_applied_1988,rivlin_chebyshev_2020,thacher1964conversion}. \begin{lemma}[\cite{fox1968chebyshev}]\label{lem:monomialstocheb} Let $n\in \mathbb{N}_0$ and $\left(\begin{smallmatrix} n\\ k \end{smallmatrix}\right)$ denote the binomial coefficient, then
\begin{align*}
x^j = \sum_{k=0}^j \gamma_{j,k} T_k(x),
\end{align*}
where
\begin{align*}
    \gamma_{j,k} = \begin{cases}
        2^{-j} \binom{j}{(j-k)/2}, \quad \hspace{1.5mm} \text{if $(j-k)$ even and $k = 0$,} \\
        2^{1-j} \binom{j}{(j-k)/2}, \quad \text{if $(j-k)$ even and $k\neq 0$,} \\
        0,  \qquad \qquad \quad \hspace{4mm} \text{else.}
    \end{cases}
\end{align*}
\end{lemma}

Using this procedure to compute Chebyshev coefficients is typically ill-advised in the context of function approximations since it involves first computing a non-Chebyshev series which can lead to loss of many of the Chebyshev polynomials' advantageous approximation properties and often requires higher precision arithmetic for sensible results, cf. \cite{cody1970survey}. However, in our method the degree of the polynomial is always fixed to the number of intersecting surfaces $m$, meaning that the required form of the Chebyshev linear combinations of the ESPs can be symbolically precomputed and directly evaluated instead of numerically computing the ESPs first.

\begin{algorithm}
    \caption{Pointwise reconstruction from smooth interpolants from the Chebyshev colleague matrix}\label{alg:chebmethod}\vspace{2mm}
    \KwIn{Chebyshev coefficients $b_j(\mathbf{x}')$ such that $$T_m(y)+\sum_{j=0}^{m-1} b_j(\mathbf{x}') T_j(y) = \kappa \sum_{k=0}^{m} (-1)^{m-k} s_{m-k}(\mathbf{f}(\mathbf{x}')) y^k$$ with arbitrary $\kappa\in \mathbb{R} \setminus \{ 0 \}$ for an underlying multi-surface $\mathbf{f}(\mathbf{x})$ with $m$ entries at a given point $\mathbf{x}' \in \Omega \subset \mathbb{R}^d$.}
    \KwOut{Pointwise values of multi-surface $\mathbf{f}(\mathbf{x}')$ sorted by value.} 
      $A^{\mathrm{C}}_p \gets$ Colleague matrix of $b_j$ via Theorem \ref{thm:colleaguedef}\\
       $ \lambda \gets \text{eig}(A^{\mathrm{C}}_p)$
        
    \Return $\text{sort}(\lambda)$
\end{algorithm}

\section{Sensitivity analysis}\label{sec:numericalanalysis}
An important result for understanding the sensitivity of our reconstruction methods to perturbations in exact arithmetic is the following classical theorem on the conditioning of polynomial root-finding in Chebyshev and monomial power bases.

\begin{theorem}\label{thm:rootconditioning}\cite[Theorem 4.1]{boyd2002computing}
    Let $r$ denote a real-valued root of multiplicity $\ell$ on the interval $x \in [-1,1]$ of a polynomial
\begin{align*}
p(x) = \sum_{j=0}^{n} b_j \phi_j(x),	\quad x \in [-1,1],
\end{align*}
where either $\phi_j(x) = T_j(x)$ or $\phi_j(x) = x^j$, both of which satisfy $|\phi_j(x)|\leq 1$ for all $x \in [-1,1]$. Furthermore, let $\tilde{r}$ denote the root of the perturbed polynomial $\tilde{p}(x)$ with modified $k$-th coefficient:
\begin{align*}
\tilde{p}(x) = (b_k + \epsilon) \phi_k(x) + \sum_{j=0, j\neq k}^n b_j \phi_j(x),	\quad x \in [-1,1].
\end{align*} 
Then, the shift of the root caused by the perturbation $0 < \epsilon < 1$ is bounded by
\begin{align*}
    |r - \tilde{r}| \leq \epsilon^{\frac{1}{\ell}} \Big| \frac{1}{\ell!} \frac{d^\ell p}{dx^\ell}(r) \Big|^{-\frac{1}{\ell}} + O\left(\epsilon^{\frac{\ell+1}{\ell}}\right).
\end{align*}
\end{theorem}

\begin{remark}
    Abstractly, without taking numerical errors during the root-finding into account, Theorem \ref{thm:rootconditioning} tells us that our method will have errors of order $\sqrt{\epsilon}$ where two surfaces cross and of order $\epsilon$ where no surfaces cross, where $\epsilon$ is the local absolute error in our interpolation of the ESPs or Chebyshev coefficients.
\end{remark}
    \begin{remark}
    Theorem \ref{thm:rootconditioning} also shows that the errors of distinct roots are not coupled, i.e. are not affected by multiplicities other than their own. This means that there is no pollution to the other surfaces at points where two surfaces cross.
\end{remark}

In the following sections we collect further results which can be used to understand the error incurred from the Schmeisser and Chebyshev colleague reconstruction approaches. The sensitivity analysis necessarily differs between these methods as the Schmeisser matrix is something we have to construct using polynomial division, cf. Algorithm \ref{alg:constructiveschmeisser}, which is an equivalent procedure to deconvolution. Since deconvolution is widely understood to be a numerically unstable process, we will focus the stability analysis for the Schmeisser companion matrix on perturbations on the companion matrices themselves, where meaningful error bounds can be given. In a practical application, however, care must be taken that sufficient precision is available for the polynomial division to proceed which may be done by using higher precision arithmetic, symbolic calculations or pre-processing steps as appropriate. In contrast, for the Chebyshev colleague matrix approach a classical result allows us to consider both perturbations to the ESP coefficients and the colleague matrix simultaneously, and thus effectively fully characterize the sensitivity to noise in the data.

\subsection{Errors in the Chebyshev colleague matrix method} \label{sec:chebanalysis}
The Chebyshev colleague matrix approach shares a key strength with the Frobenius approach in that once the coefficients $b_j(\mathbf{x})$ discussed in Theorem \ref{thm:colleaguedef} are computed there is no further error accumulation in the construction of the colleague matrix. A major drawback compared to the Schmeisser approach, however, is that the colleague matrix is in general not symmetric and thus the algorithm cannot guarantee real eigenvalues being returned for sufficiently large perturbations. To understand the sensitivity of the Chebyshev colleague reconstruction method we can leverage results for the error control of Chebyshev series root-finding.

\begin{theorem} \cite[Corollary 5.4]{noferini_structured_2021}\label{thm:chebcollperturbation}
    Let $A^{\mathrm{C}}_p = H - e_1 c^\top$ with $H = H^\top$ be the Chebyshev colleague matrix of a polynomial $p(x) = T_n(x) + \sum_{j=0}^{n-1} b_j T_j(x)$ with  $b_j \in \mathbb{R}$ as in Theorem \ref{thm:colleaguedef} and $||\delta H ||_2 \leq \epsilon_H$, $||\delta e_1|| \leq \epsilon_1$ and $||\delta c|| \leq \epsilon_c$. Then the perturbed matrix
    \begin{align*}
        A^{\mathrm{C}}_p + \delta A^{\mathrm{C}}_p = H + \delta H - \frac{1}{2} (e_1 + \delta e_1) (c + \delta c)^\top
    \end{align*}
    is the colleague matrix corresponding to a polynomial
    \begin{align*}
        p(x) + \delta p(x) = \sum_{j=0}^n (b_j + \delta b_j) T_j(x),
    \end{align*}
    with perturbed coefficients $b_j + \delta b_j$, for which we have the bounds
    \begin{align*}
       |\delta b_j| \leq (6 \|c\|_2 \epsilon_1 + 2 \sqrt{n} \epsilon_c + (5+16\sqrt{n}\|c\|_2)\epsilon_H )n^2 + \mathcal{O}(\epsilon_H^2+\epsilon_1^2+\epsilon_c^2).
    \end{align*}
\end{theorem}
\begin{remark}
    In the context of our problem, perturbations occur at the level of the data and thus the ESPs to be interpolated. The perturbations $\epsilon_1$ and $\epsilon_H$ are of the order of numerical precision and can normally be neglected. The bound from Theorem \ref{thm:chebcollperturbation} then simply reads
    \begin{align*}
       |\delta b_j| \leq 2 n^{5/2} \| \delta c\|  + C \| \delta c \|^2,
    \end{align*}
    where $C$ may depend on $\|c\|$ and $n$.
\end{remark}

If a backward stable algorithm is used to compute the eigenvalues of the colleague matrix, then the root-finding is backward stable as a function of perturbations of the Chebyshev series' coefficients, cf. \cite[Cor. 2.8, Rem. 2.9]{noferini2017chebyshev} and \cite{serkh_provably_2021}. Since no deconvolution is needed for the construction of the Chebyshev colleague matrix, the approach based on it is expected to perform well in generic circumstances.

\subsection{Errors in the Schmeisser companion matrix method}
Error bounds for the numerical eigenvalue problem associated with symmetric tridiagonal matrices have been a topic of long standing research~\cite{golub1962bounds,wilkinson1958calculation}. More recent results~\cite{hladik2017eigenvalues,hladik2011characterizing, hladik2010bounds} allow for the characterization of upper and lower bounds for the eigenvalues of symmetric tridiagonal interval matrices which are defined as follows.

\begin{definition}[Symmetric tridiagonal interval matrix]
We define the symmetric tridiagonal interval matrix $\mathbf{A} = [\underline{A}, \overline{A}]$ to be the set of all symmetric tridiagonal matrices $A_j$ with entries taken from two given intervals, that is with all $a_j \in [\underline{a}_j,\overline{a}_j]$ with $j = 1,...,n$ and all $b_j \in [\underline{b}_j,\overline{b}_j]$ with $j=2,...,n$.
\end{definition}

We can think of interval matrices as representing a matrix with error bounds on each of its entries. These objects are studied in the field of interval arithmetic, providing an elegant way for rigorous numerics with guaranteed error bounds~\cite{moore2009introduction,tucker2011validated,mayer2017interval}. We seek error bounds on the eigenvalues, i.e. eigenvalue intervals $\bm{\lambda}_j=[\underline{\lambda}_j,\overline{\lambda}_j]$, of an interval matrix $\mathbf{A}$. Note that any method to determine the upper bounds $\overline{\lambda}_j$ could be used to determine the lower bounds via the transformation $A \mapsto -A$, cf. \cite{hladik2017eigenvalues,hladik2010bounds}. The following result, stated without proof in \cite{rohn2005handbook} and proved in \cite{hladik2010bounds}, establishes a key result.

\begin{theorem}\label{thm:symtriintervalarithmetic} \cite[Theorem 3.1]{hladik2010bounds}
Let $\bm{\lambda}_j(\mathbf{A})$ denote the $j$-th eigenvalue interval of the symmetric tridiagonal interval matrix $\mathbf{A} = [\underline{A},\overline{A}]$.
We denote by $A_c$ and $A_\Delta$ the mid point and radius matrices of the interval matrix defined by
\begin{align*}
A_c := \frac{1}{2}\left(\overline{A} + \underline{A} \right),\qquad
A_\Delta := \frac{1}{2}\left( \overline{A} - \underline{A}\right).
\end{align*}
Then
\begin{align*}
    \bm{\lambda}_j(\mathbf{A}) \subseteq [\lambda_j(A_c)-\rho(A_\Delta),\lambda_j(A_c)+\rho(A_\Delta)],
\end{align*}
where $\rho(A_\Delta) := \max\{|\lambda_1(A_\Delta)|,...,|\lambda_n(A_\Delta)|\}$ denotes the spectral radius of $A_\Delta$.
\end{theorem}

In the context of the Schmeisser method discussed in Section \ref{sec:method}, the natural interpretation of the interval matrix $\mathbf{A} = [\underline{A},\overline{A}]$ is as $$\mathbf{A} = A \pm \mathcal{E}_{c} = [A - \mathcal{E}_{c}, A +  \mathcal{E}_{c}],$$
where $\mathcal{E}_{c}$ is an upper bound for the absolute error accumulated in the construction of the Schmeisser companion matrix via Algorithm \ref{alg:constructiveschmeisser}. Defined in this way, the radius matrix $A_\Delta = \mathcal{E}_{c}$ is simply an absolute error bound matrix. We will consider an upper bound over the error terms for each of the bands considered in $\mathcal{E}_{c}$, which has the notable advantage that its eigenvalues and thus also its spectral radius are explicitly known by the following results.
\begin{lemma}\label{lem:classicaltoeplitzeigvals} [\cite[p.86 and p.113]{smith1985numerical} or \cite{noschese2013tridiagonal}]
    Let $T_n$ be a general $n \times n$ tridiagonal Toeplitz matrix, i.e.
    \begin{align*}
        T_n = \begin{pmatrix}
        a && b && 0 && \hdots && 0\\
        c && a && b && \ddots && 0\\
        0 && c && a && \ddots && 0\\
        \vdots && \ddots && \ddots && \ddots && b\\
        0 && \hdots && 0 && c && a
    \end{pmatrix},
    \end{align*}
    where $a, b, c$ may be real or complex. Then the eigenvalues of $T_n$ are explicitly given by
    \begin{align*}
        \lambda(T_n)_j &= a + 2 \sqrt{bc} \cos\left( \tfrac{j \pi}{n+1} \right) = a + 2 \sqrt{|bc|} e^{i\frac{\mathrm{arg}(c)+\mathrm{arg}(b)}{2}} \cos\left( \tfrac{j \pi}{n+1} \right).
    \end{align*}
\end{lemma}
\begin{corollary}\label{cor:spectralradiusexplicit}
Let $T_n$ be a tridiagonal Toeplitz matrix as in Lemma \ref{lem:classicaltoeplitzeigvals}. Then the spectral radius of $T_n$ satisfies
\begin{align*}
    \rho(T_n) &= \max\left\{\Big| a + 2 \sqrt{bc} \cos\left( \tfrac{\pi}{n+1} \right) \Big|,\Big|a + 2 \sqrt{bc} \cos\left( \tfrac{n \pi}{n+1} \right)\Big|\right\} \leq | a | + 2 | \sqrt{bc} |.
 \end{align*}
If $T_n$ is real symmetric, i.e. $c=b$ and $a,b\in\mathbb{R}$, then this further simplifies to 
\begin{align*}
    \rho(T_n) &= \max\left\{\Big| a + 2 |b| \cos\left( \tfrac{\pi}{n+1} \right) \Big|,\Big|a + 2 |b| \cos\left( \tfrac{n \pi}{n+1} \right)\Big|\right\} \leq | a | + 2 |b |.
 \end{align*}
\begin{proof}
    The explicit expressions are an immediate corollary of Lemma \ref{lem:classicaltoeplitzeigvals} and the definition of the spectral radius $\rho(T_n) := \max\{|\lambda_1(T_n)|,...,|\lambda_n(T_n)|\}$, cf. \cite{noschese2013tridiagonal}. The upper bounds follow from the subadditivity and multiplicativity of the absolute value function and $\forall \theta \in \mathbb{R}:|\cos(\theta)|\leq 1$.
\end{proof}
\end{corollary}
\begin{remark}
    The explicit spectral radius for a tridiagonal Toeplitz matrix $T_n$ given in Corollary \ref{cor:spectralradiusexplicit}  depends on the size of the $n \times n$ matrix, while the bounds do not. By monotonicity and since
    \begin{align*}
       \lim_{n \rightarrow \infty}\cos\left( \tfrac{\pi}{n+1} \right) = 1,\qquad
       \lim_{n \rightarrow \infty}\cos\left( \tfrac{n \pi}{n+1} \right) =  -1.
    \end{align*}
    one observes that the bounds become sharper with increasing $n$ for fixed $a,b$.
\end{remark}

With the above results established, we can state the full error bound for computing the eigenvalues of a Schmeisser companion matrix.

\begin{theorem}\label{thm:mainerrorthm}
        Let $A(\mathbf{x})$ be a symmetric tridiagonal Schmeisser companion matrix with real, non-negative entries at each point $\mathbf{x} \in \Omega \subset \mathbb{R}$  and let $\mathcal{E}_c$ be a symmetric tridiagonal Toeplitz matrix with $a=b=c= \varepsilon_c \geq 0$. Then for any family of matrices $\hat{A}(\mathbf{x}) \in [A(\mathbf{x})-\mathcal{E}_c,A(\mathbf{x})+\mathcal{E}_c]$ we have       
                \begin{align*}
            \lambda_j(\hat{A}(\mathbf{x})) \in \bm{\lambda}_j([A(\mathbf{x})-\mathcal{E}_c,A(\mathbf{x})+\mathcal{E}_c]) \subseteq [\lambda_j(A(\mathbf{x}))-3\varepsilon_c,\lambda_j(A(\mathbf{x}))+3\varepsilon_c]
        \end{align*}
        and furthermore
        \begin{align*}
            \sup_{\mathbf{x} \in \Omega}|\lambda_j(A(\mathbf{x}))-\lambda_j(\hat{A}(\mathbf{x}))| \leq 3 |\varepsilon_c|.
        \end{align*}
\begin{proof}
    The statements follow by combining Corollary \ref{cor:spectralradiusexplicit} with Theorem \ref{thm:symtriintervalarithmetic}. 
\end{proof}
\end{theorem}

Theorem \ref{thm:mainerrorthm} establishes that the method in Algorithm \ref{alg:schmeissermethod} converges uniformly as a function of perturbations on the elements of the Schmeisser companion matrix. It is important to remind ourselves that this unfortunately does not necessarily translate to well-behaved convergence in terms of the error in an approximation of the ESPs and underlying data, since the construction of the Schmeisser companion matrix involves the numerically ill-conditioned step of polynomial division, which is mathematically equivalent to a deconvolution.

There are thus two primary drawbacks to the Schmeisser approach: First, the construction algorithm will fail if it is started with a ESP polynomial perturbed in such a way that its roots are no longer strictly real (though this can be substantially alleviated) and second, while there are tight error bounds for the eigenvalue computation once the Schmeisser matrix is constructed, the stability of the Schmeisser construction algorithm itself is in doubt due to the equivalence of polynomial division and deconvolution. It remains to be seen whether the latter shortcoming can be addressed without costly use of higher precision arithmetic in the construction process, e.g. by pre-processing steps.

\section{Numerical experiments}\label{sec:experiments}
In this section we illustrate the performance of the proposed reconstruction methods in various toy problems designed to showcase generic interesting behavior with respect to convergence and stability as well as two applications-oriented examples with more complex behavior. In the context of the relevant applications, the most widely used error concepts are those of mean absolute error (MAE) and root mean square error (RMSE) rather than the more strict maximum absolute error often encountered in the context of polynomial approximations. While very useful for discussing worst-case behavior in toy problems with full control of the data, maximum absolute error is much more sensitive to outliers and thus less useful when investigating the different error behaviors in the application-adjacent examples. On the other hand, MAE and RMSE risk masking the extent to which the direct reconstruction method performs poorly by averaging out the error over the rest of the domain.

Thus, while we still also present these more conventional errors, we will introduce a gap-weighted error to discuss the specific behavior of interest which is the resolution of surface crossings in a multi-surface. Denoting by $$\Delta_{ij} \mathbf{f}(\mathbf{x}) := \mathbf{f}(\mathbf{x})_j - \mathbf{f}(\mathbf{x})_i,$$
the difference between the $i$-th and $j$-th components of a multisurface $\mathbf{f}$ at $\mathbf{x}$, we introduce the following gap-weighted error for approximate reconstructions $\tilde{\mathbf{f}}$ \begin{align} \label{eq:gapweightederror}
    \mathrm{err}_{\Delta}(\mathbf{f},\tilde{\mathbf{f}}) = \max_{x \in \Omega, i\neq j} \frac{| \Delta_{ij}\tilde{\mathbf{f}}(\mathbf{x})-\Delta_{ij}\mathbf{f}(\mathbf{x}) |}{\epsilon_W+ |\Delta_{ij} {\mathbf{f}}(\mathbf{x})|}
\end{align}
where $\epsilon_W>0$ is a suitably small parameter to prevent the denominator from reaching $0$ at crossings which we interpret as a target accuracy. This error places additional weight on good approximations near cusps as $\Delta_{ij} {\mathbf{f}}(\mathbf{x}) \rightarrow 0$ as $\mathbf{x}$ approaches a crossing site. We note that ${\rm err}_\Delta$ is invariant under permutations of the surfaces. For all our numerical experiments we choose $\epsilon_W = 5 \times 10^{-2}$. 

\subsection{Toy problems}
\subsubsection{Smooth intersecting sinusoid surfaces}\label{sec:numericaltoy1}
\begin{figure}\centering
     \subfloat[]
    {{\includegraphics[width=7.3cm]{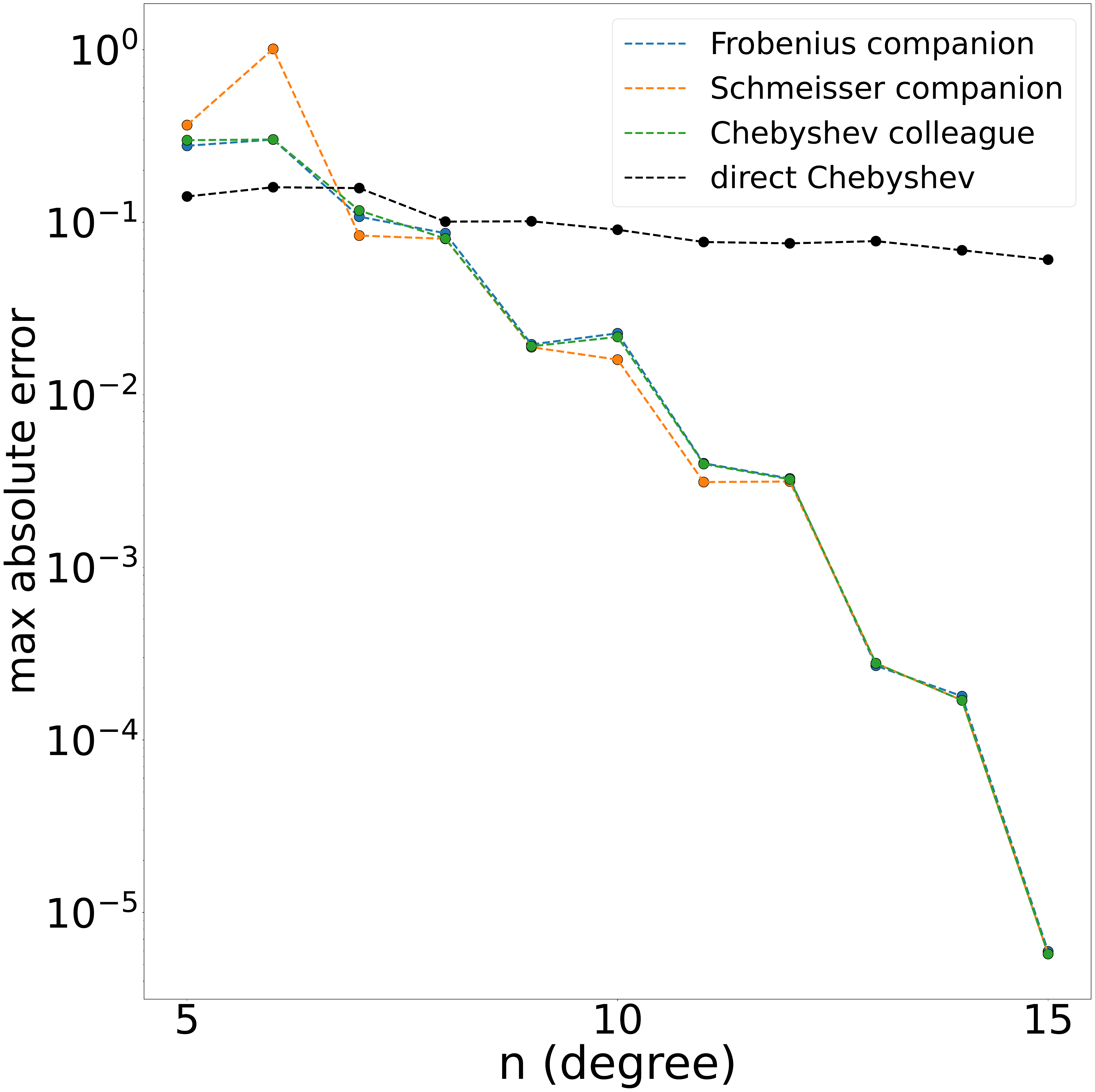}}}
         \subfloat[]
    {{\includegraphics[width=7.3cm]{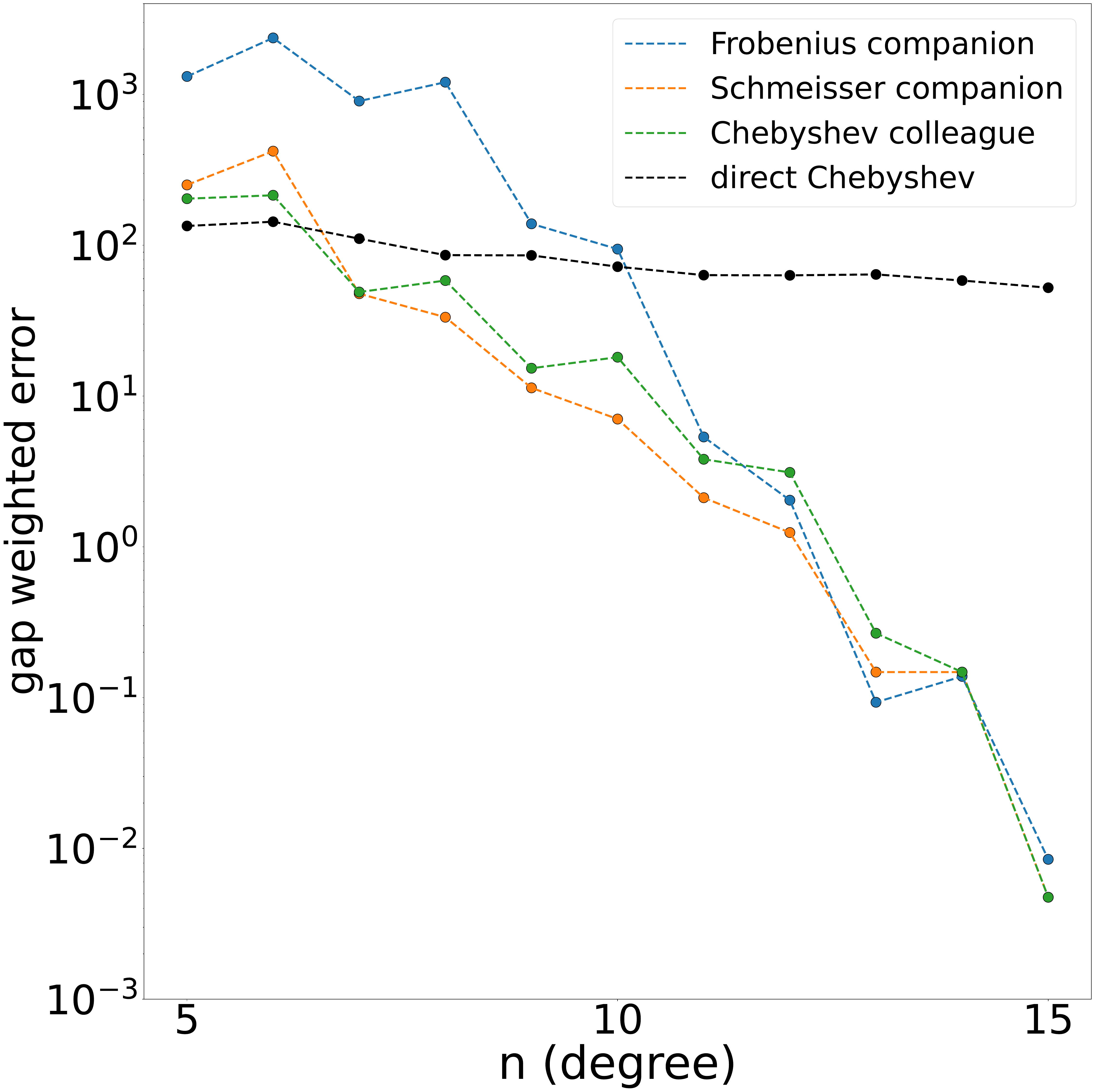}}}\\
     \subfloat[]
    {{\includegraphics[width=7.3cm]{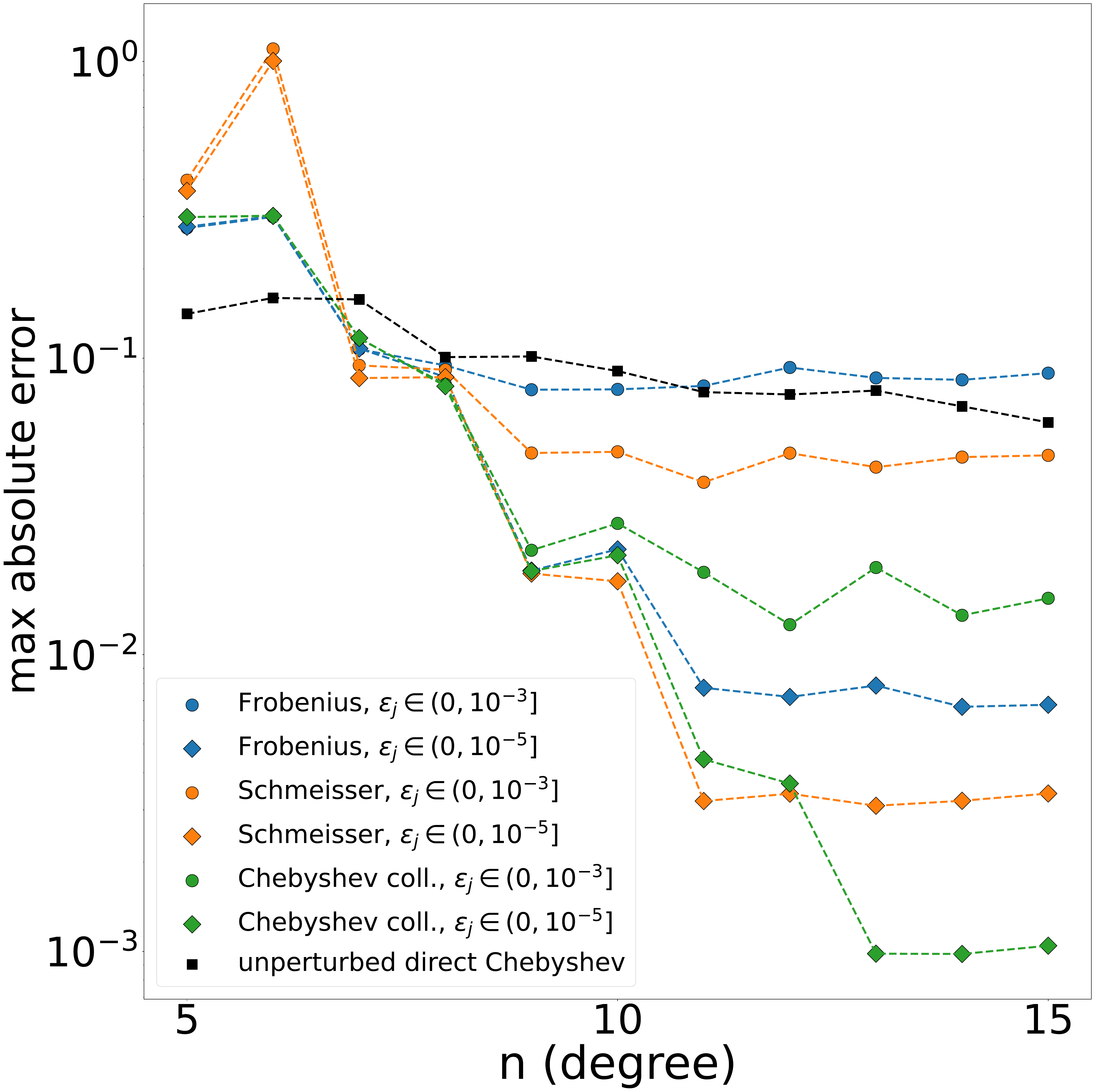}}}
         \subfloat[]
    {{\includegraphics[width=7.3cm]{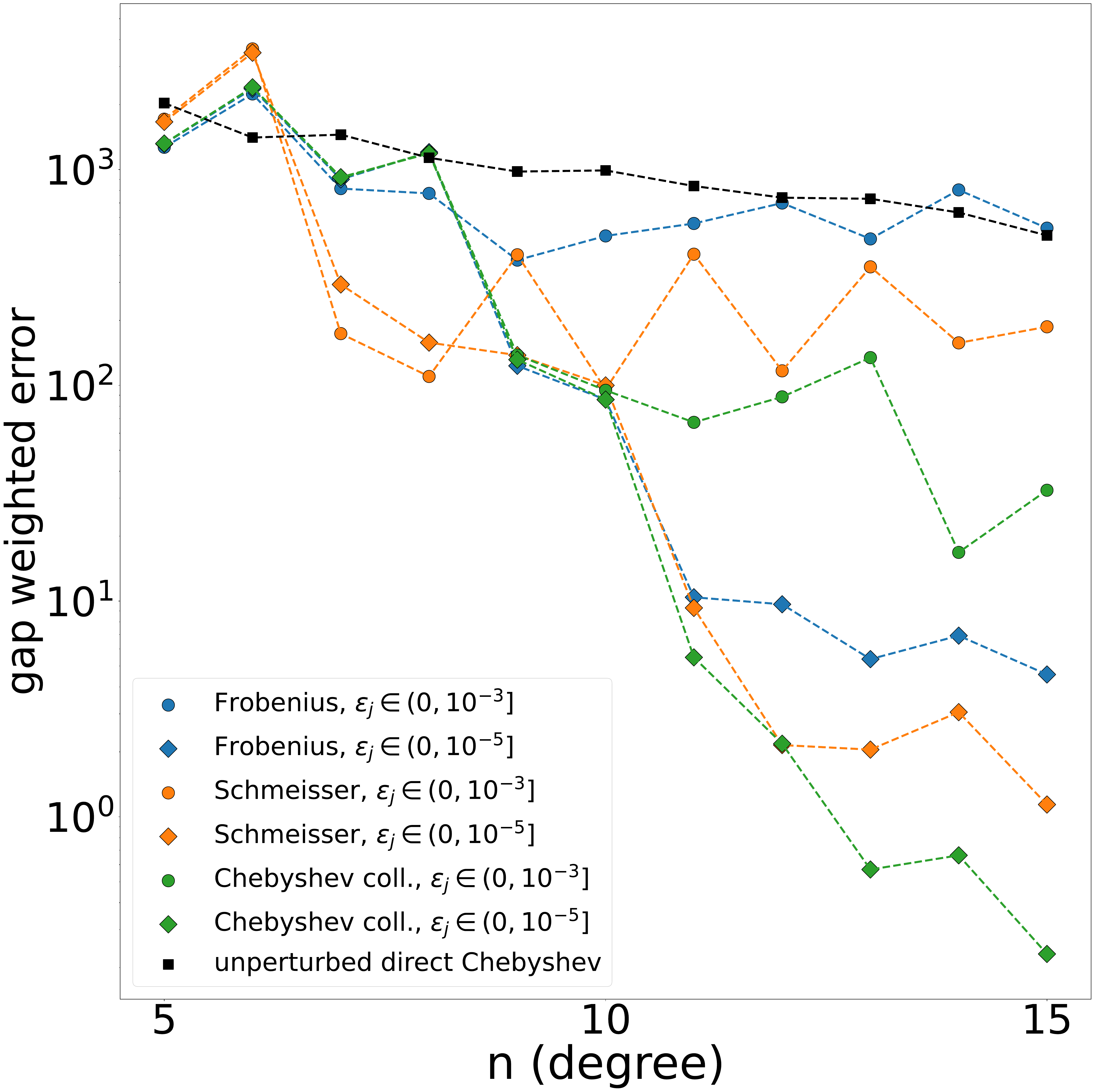}}}
    \caption{(a-b) Semi-logarithmic maximum absolute error and gap weighted error on all $1000$ random gridpoints (sampled uniformly) for the problem in Eqs. (\ref{eq:toy1start}--\ref{eq:toy1end}) solved via Frobenius, Schmeisser companion, Chebyshev colleague or direct Chebyshev interpolation approaches. (c-d) show analogous error plots with added perturbations of indicated magnitude on the surface or ESP evaluations.}
    \label{fig:1dtoyerrors}
\end{figure}

\begin{figure}\centering
\includegraphics[width=6.3cm]{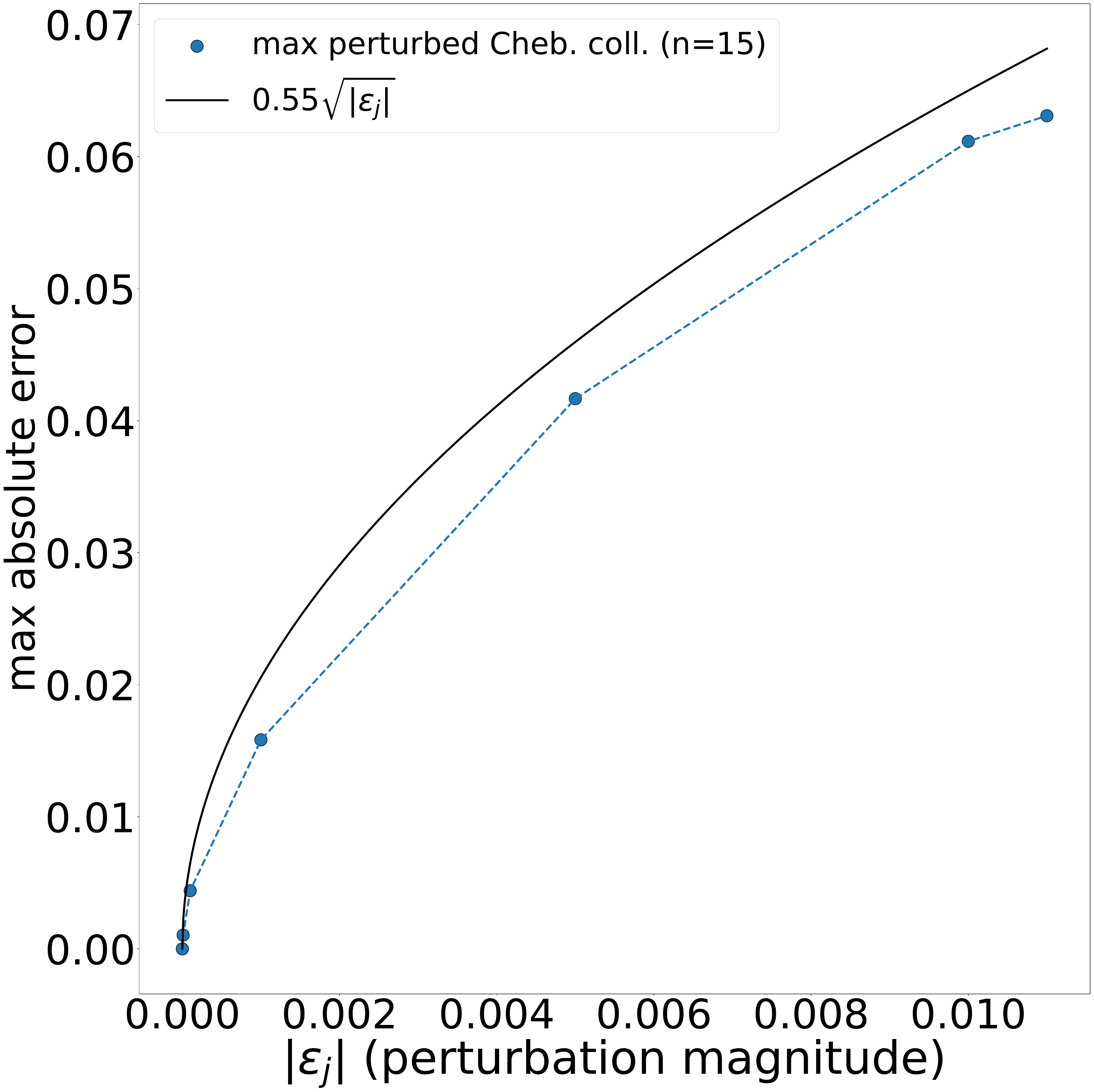}
    \caption{Endpoint plot of Figure \ref{fig:1dtoyerrors}(c) including some additional data points, showing that the converged (in degree) error for using the Chebyshev colleague approach for the problem in (\ref{eq:toy1start}--\ref{eq:toy1end}) scales like the square root of the perturbations, consistent with Theorem \ref{thm:rootconditioning} in light of the multiplicity $2$ roots observed in this problem, cf. Figure \ref{fig:firsttoyexample}.}
    \label{fig:1dtoyrooterrors}
\end{figure}

We begin with the problem shown in the introduction in Figure \ref{fig:firsttoyexample} featuring the following three intersecting sinusoid curves on $[0,2]$:
\begin{align}\label{eq:toy1start}
    z_1(x) &= \sin(x),\\ 
    z_2(x) &= \cos(2x),\\
    z_3(x) &= \sin(2x).\label{eq:toy1end}
\end{align}
Figure \ref{fig:1dtoyerrors}(a) shows maximum absolute error plots for the reconstruction using our proposed methods as well as a direct Chebyshev interpolation of the cusps as in Figure \ref{fig:firsttoyexample}(b) as a function of the interpolation degree using $1000$ randomly distributed points as data (sampled from a uniform distribution), while Figure \ref{fig:1dtoyerrors}(b) shows the gap weighted error defined in \eqref{eq:gapweightederror}. As expected from the description given in Section \ref{sec:numericalanalysis}, the companion matrix approaches achieve spectral convergence in degree given data saturation while attempting to simply resolve the cusps directly with Chebyshev polynomials results in approximately $O(n^{-1})$ convergence. Figure \ref{fig:1dtoyerrors}(c) and (d) show analogous error plots but at each point the evaluation of the ESPs (and Chebyshev equivalents) are perturbed by the indicated different magnitude random amounts $\epsilon_j$. We observe that in this more realistic scenario for machine learning applications, while the convergence rate is capped at a fixed amount, our methods nevertheless perform better once $\epsilon_j \lessapprox 0.001$ for the Schmeisser variant and $\epsilon_j \lessapprox 0.01$ for the Chebyshev colleague approach. Figure \ref{fig:1dtoyrooterrors} shows that the observed error behavior with respect to the magnitude of perturbations using the Chebyshev colleague method agrees with the prediction of Theorem \ref{thm:rootconditioning}.

As we are considering a toy problem with controlled uniform distribution perturbations without outliers or sampling bias related sources of error, the gap weighted error loosely retains the hierarchy seen in the max abs errors.

\begin{figure}\centering
     \subfloat[]
    {{ \centering \includegraphics[width=7.8cm]{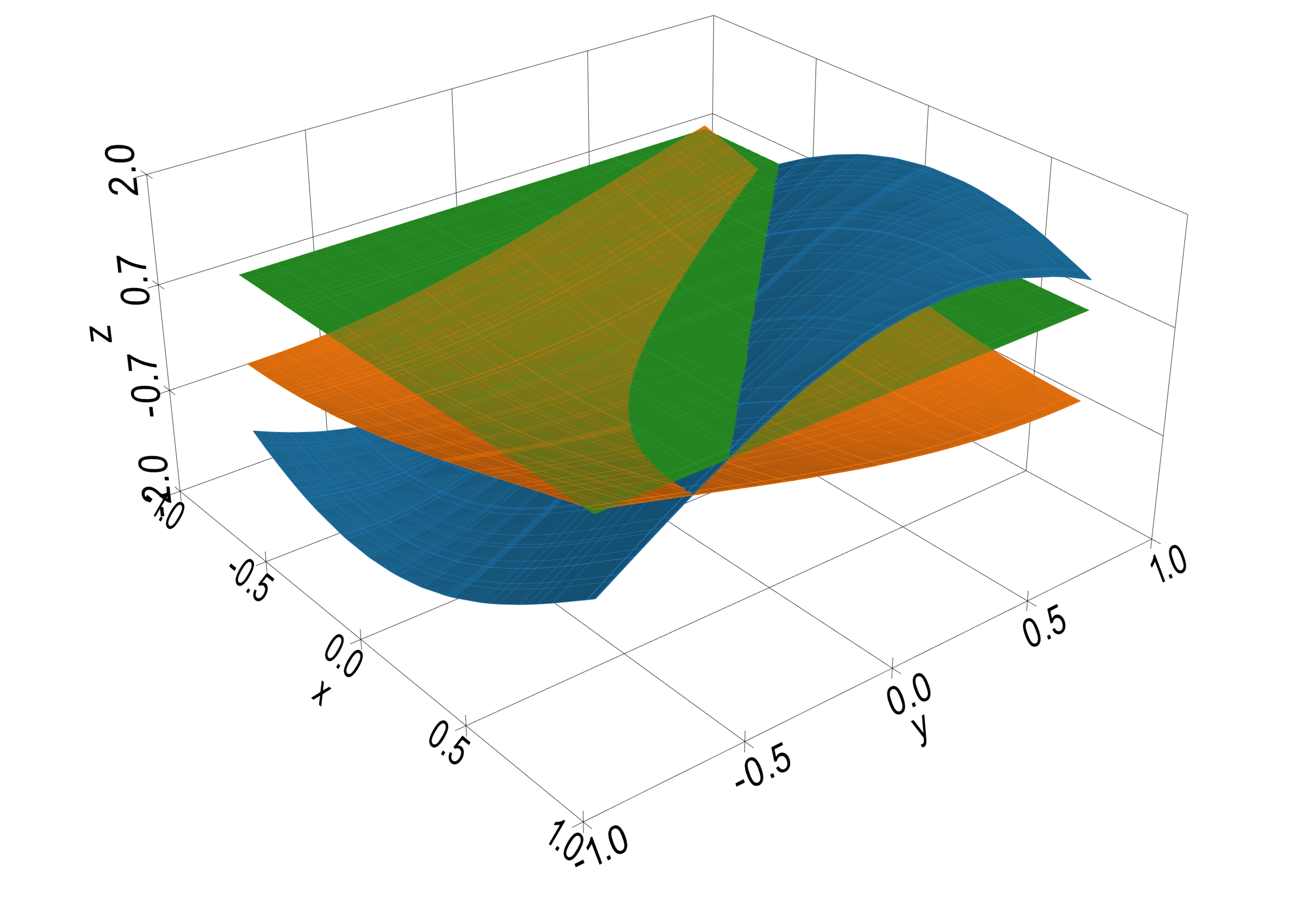} }}
     \subfloat[]
    {{ \centering \includegraphics[width=7.3cm]{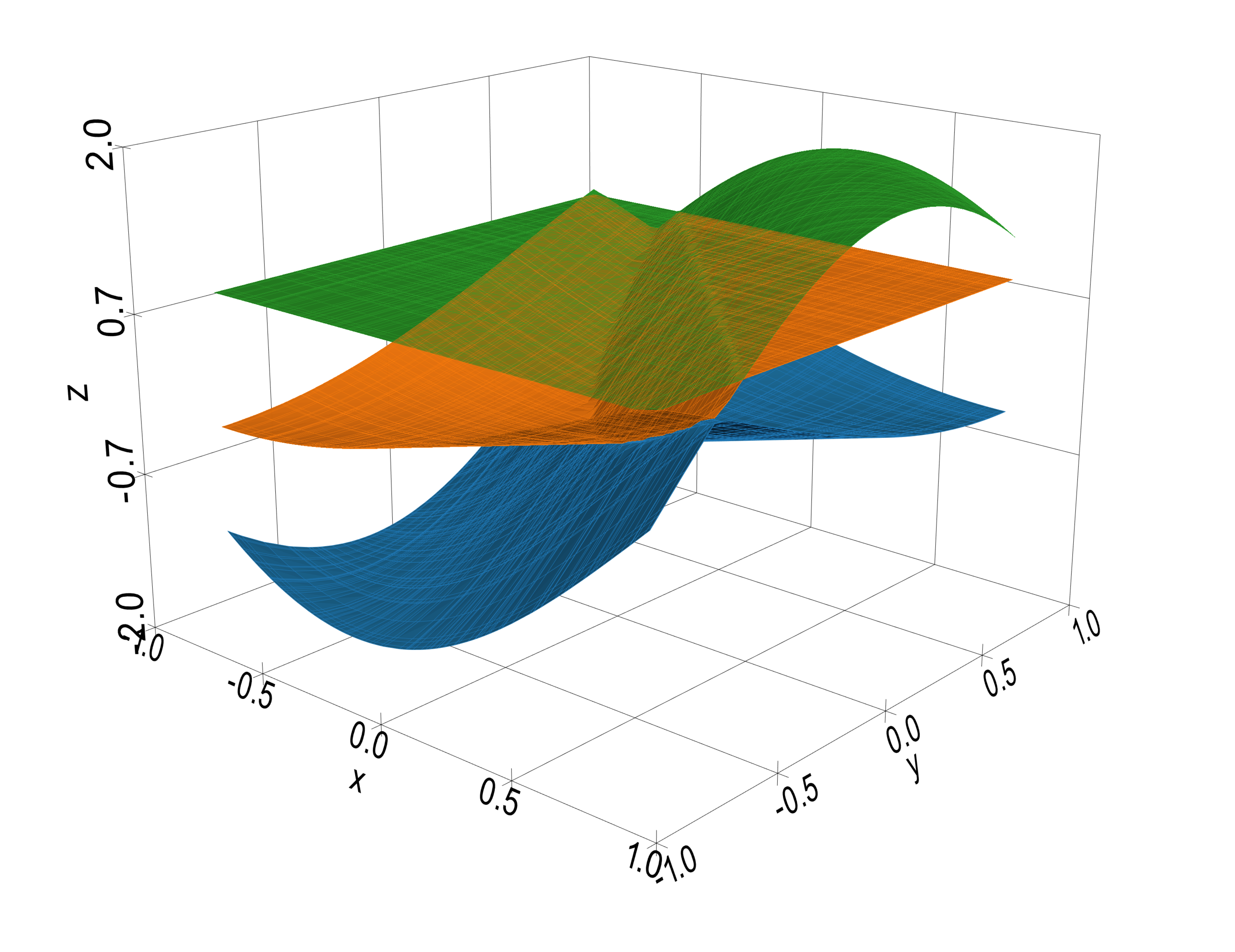} }}
    \caption{(a) shows the surfaces $\bar{\mathbf{z}}$ in (\ref{eq:2dsinusoidaltoyeq}--\ref{eq:2dsinusoidaltoyeqend}), (b) shows a corresponding value-ordered reconstruction with cusps.}
    \label{fig:2dsinusoidaltoyfig}
\end{figure}

\begin{figure}\centering
     \subfloat[]
    {{ \centering \includegraphics[width=7.3cm]{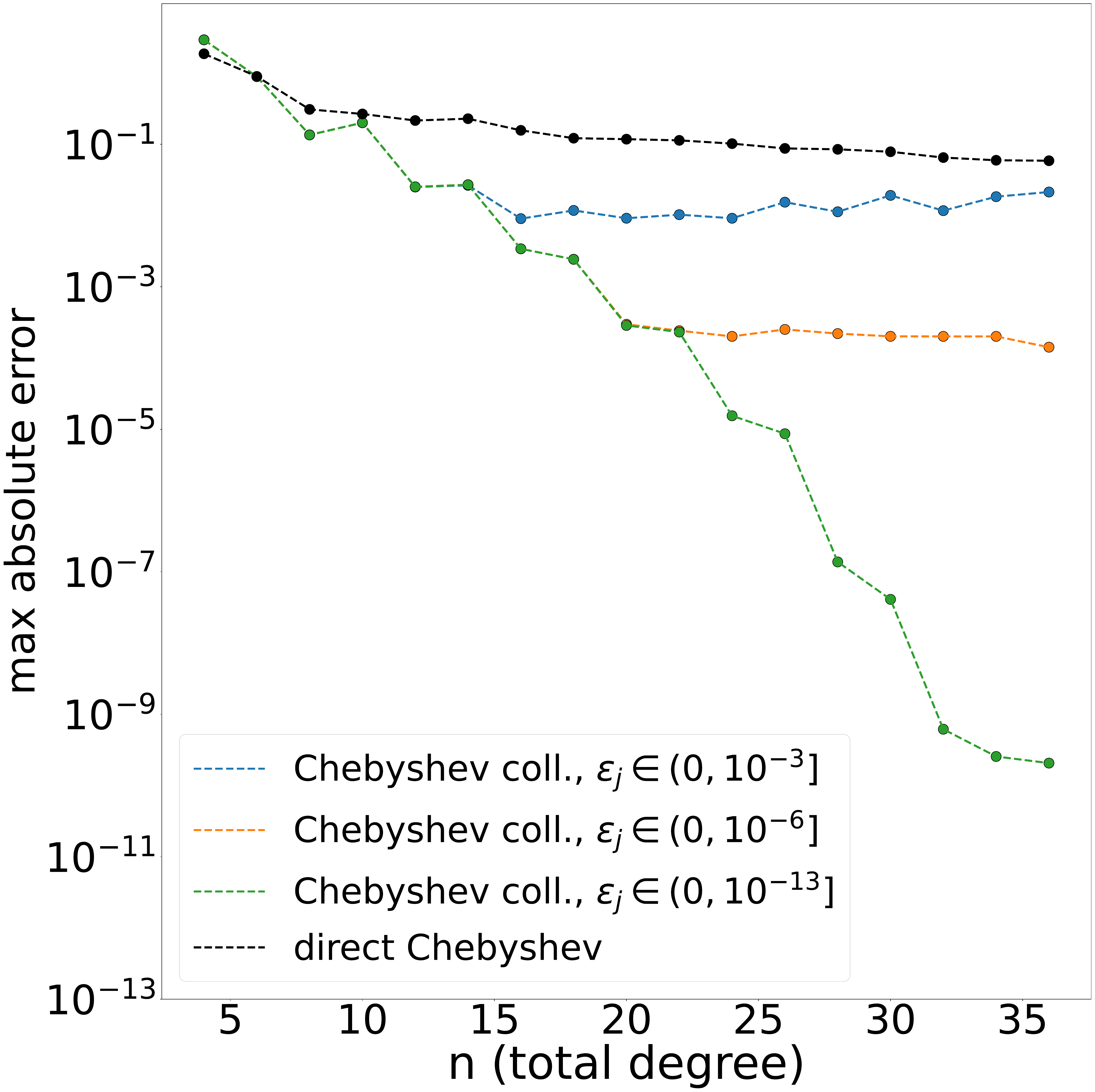} }}
    \subfloat[]
    {{ \centering \includegraphics[width=7.3cm]{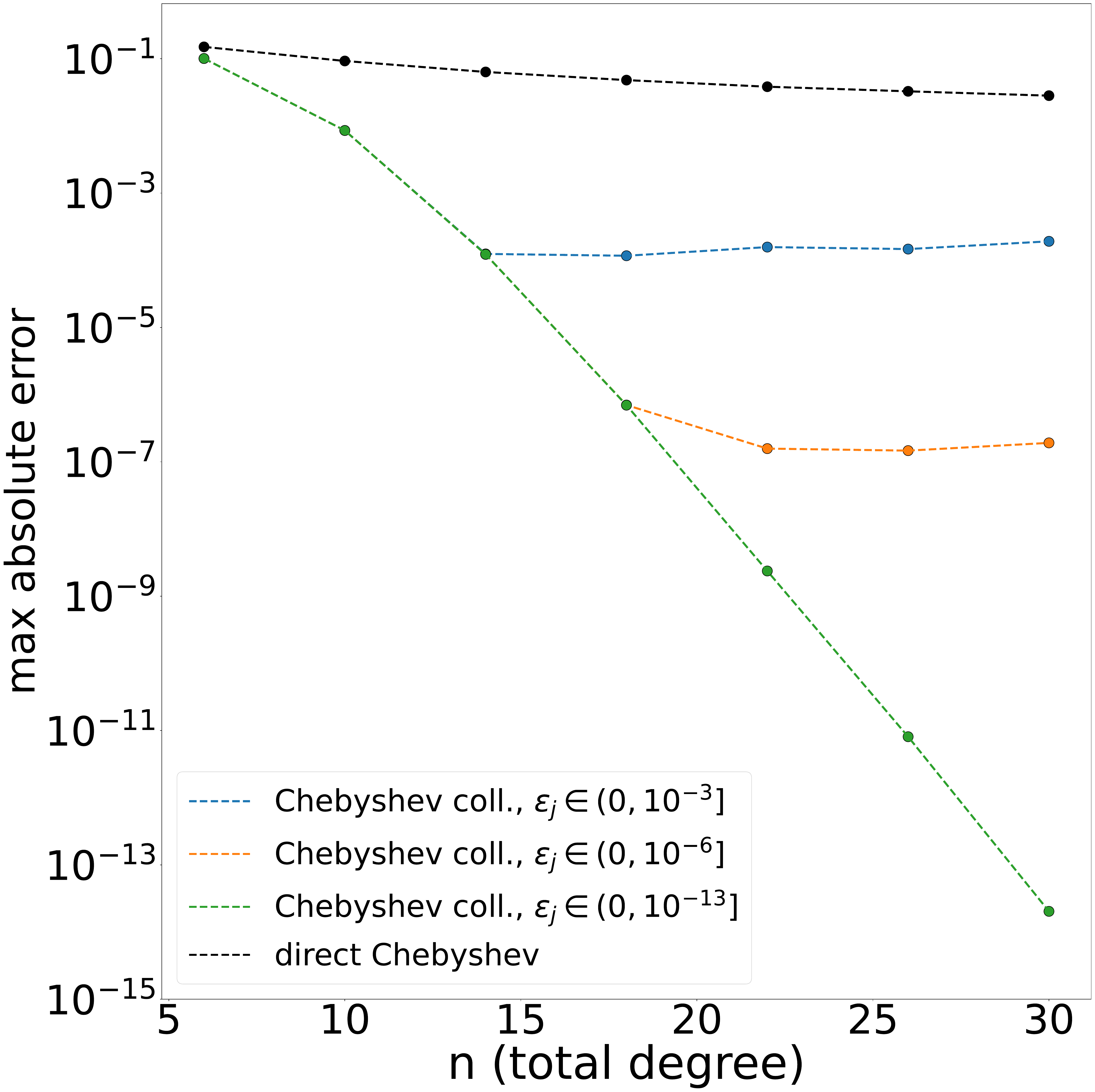} }}
    \caption{(a) shows a semi-logarithmic plot of the max. abs. error for the problem in (\ref{eq:2dsinusoidaltoyeq}--\ref{eq:2dsinusoidaltoyeqend}) using the Chebyshev colleague approach, (b) shows an analogous plot for reconstructing the conical cusp problem in (\ref{eq:toyconicalsinh}--\ref{eq:toyconicalsinhend}).}
    \label{fig:2dsinusoidaltoyfigreconstruction}
\end{figure}

Next we consider the following three intersecting 2D surfaces (cf. Figure~\ref{fig:2dsinusoidaltoyfig}):
\begin{align}\label{eq:2dsinusoidaltoyeq}
    \bar{z}_1(x,y) &= 2 \sin \left(\tfrac{6 (x+y)}{5}\right),\\
    \bar{z}_2(x,y) &= \tfrac{2}{3}-\cos (x-y),\\
    \bar{z}_3(x,y) &= 1. \label{eq:2dsinusoidaltoyeqend}
\end{align}
We plot the max. abs. error of the direct reconstruction, Chebyshev colleague and Schmeisser companion approaches, including various levels of perturbation, in Figure~\ref{fig:2dsinusoidaltoyfigreconstruction}. 

For the sake of presentation clarity we will henceforth mostly restrict our error plots to comparisons of Chebyshev colleague matrix and direct method reconstructions as the colleague matrix approach consistently performs either better or equivalently to Schmeisser or Frobenius companion matrix approaches without additional processing. Furthermore, we will only show max. abs. errors unless a meaningful qualitative difference in behavior could be observed from showing other types of errors.
\subsubsection{Conical cusps in two dimensional surfaces}\label{sec:numericaltoy2cusp}
\begin{figure}\centering
     \subfloat[]
    {{ \centering \includegraphics[width=7.3cm]{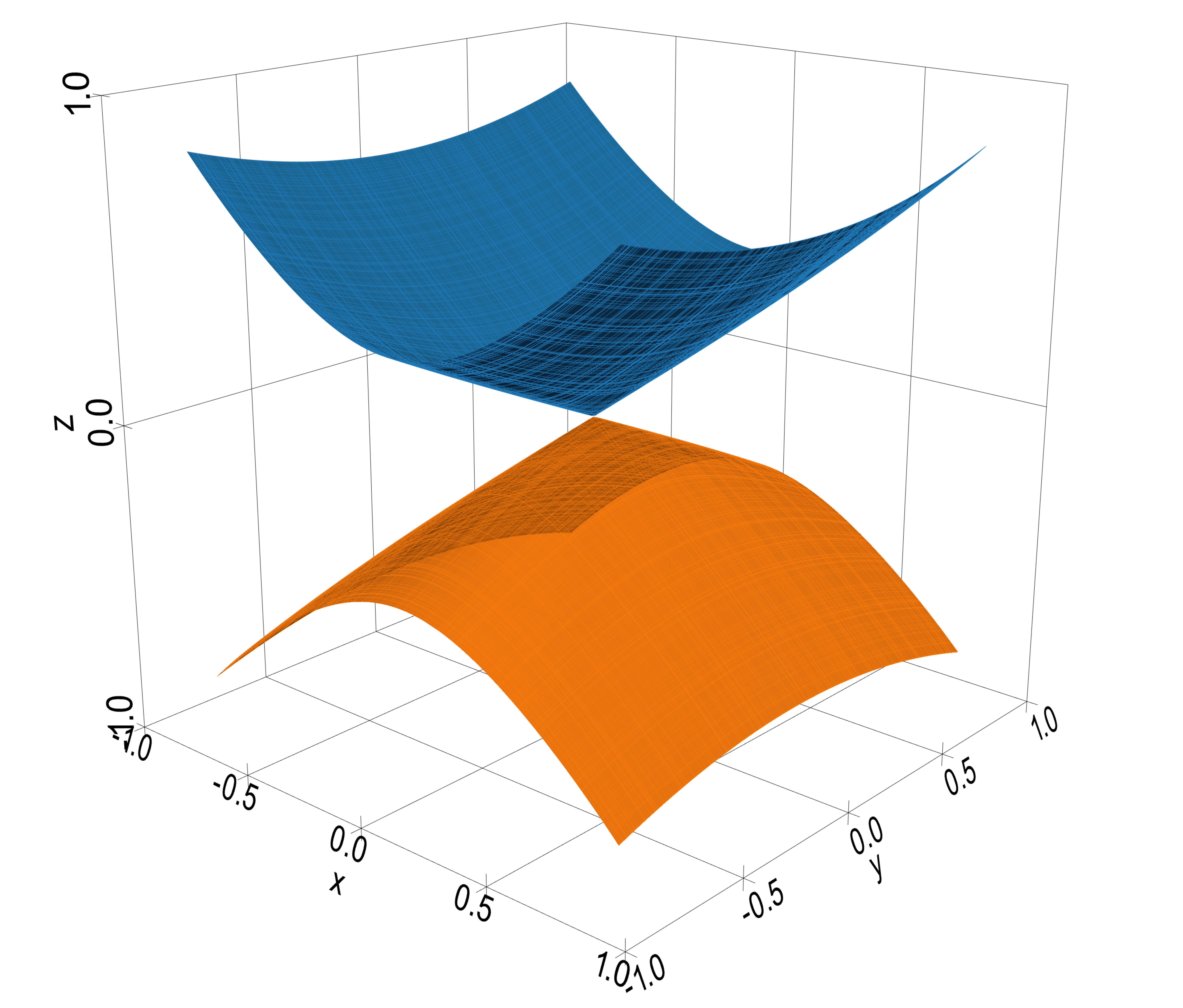} }}
     \subfloat[]
    {{ \centering \includegraphics[width=7.3cm]{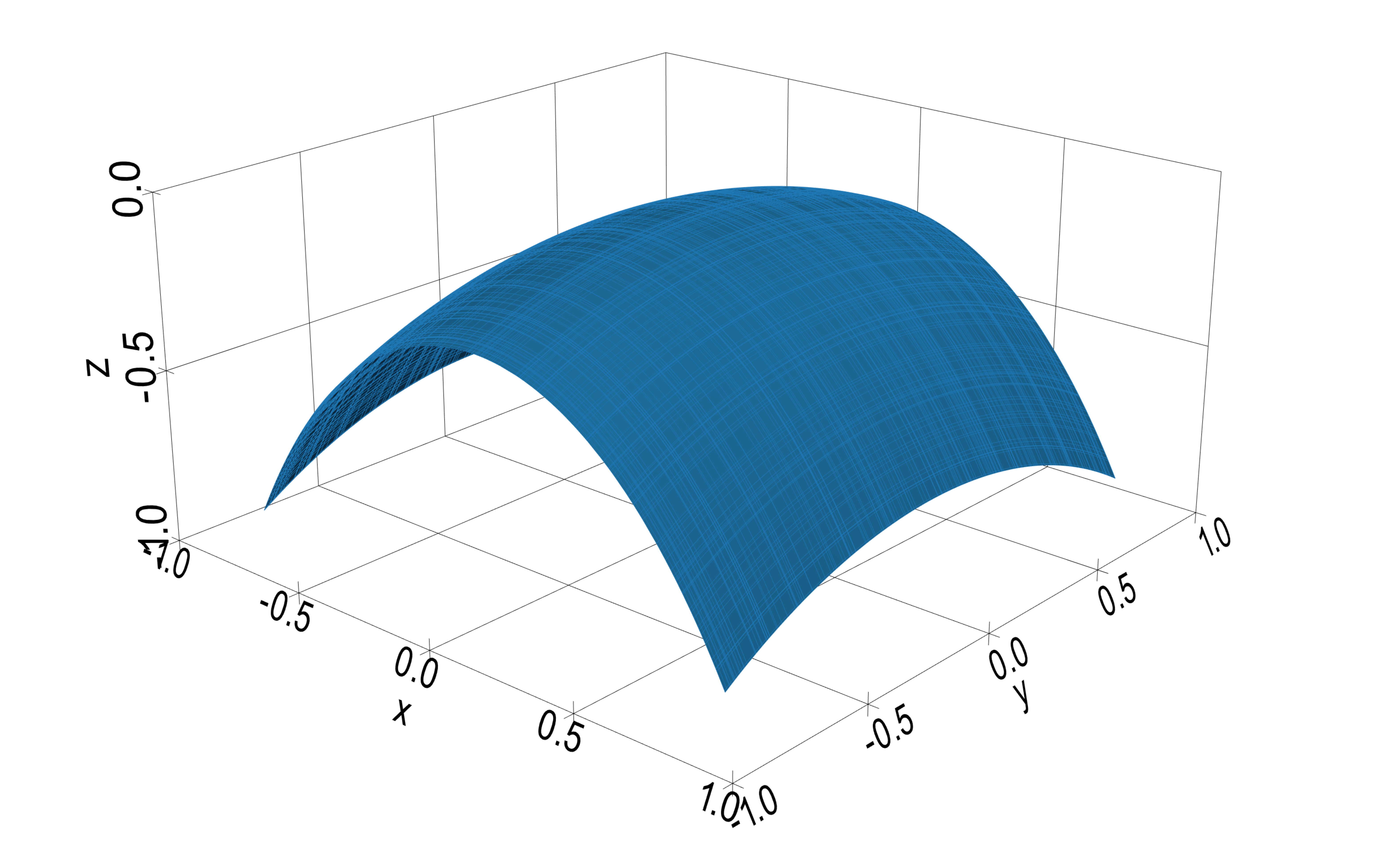} }}
    \caption{(a) shows $\bar{\mathbf{z}} = (\bar{z}_1(x,y),\bar{z}_2(x,y)$ in Eq. (\ref{eq:toyconicalsinh}--\ref{eq:toyconicalsinhend}) for $(a,b) = (\frac{4}{3}, \frac{12}{5})$, (b) shows the corresponding smooth ESP $\bar{s}_2(\mathbf{z}(x,y))$.}
    \label{fig:sinhconicaltoyexampledemo}
\end{figure}

Let two 2D surfaces $z_1$ and $z_2$ be given as the two solutions of 
\begin{align*}
    z(x,y)^2 = \left( \frac{x}{a} \right)^2 +  \left( \frac{y}{b} \right)^2,
\end{align*}
with $0 \neq a,b \in \mathbb{R}$, i.e.
\begin{align}\label{eq:toyconical}
    z_1(x,y) &= \tfrac{\sqrt{a^2 y^2+b^2 x^2}}{a b},\\
    z_2(x,y) &= -\tfrac{\sqrt{a^2 y^2+b^2 x^2}}{a b}.\label{eq:toyconicalend}
\end{align}
While neither $z_1(x,y)$ nor $z_2(x,y)$ are even once continuously differentiable, it is easy to observe that their elementary symmetric polynomials are smooth:
\begin{align*}
    s_1(\mathbf{z}(x,y)) &= z_1(x,y) + z_2(x,y) = 0,\\
    s_2(\mathbf{z}(x,y)) &= z_1(x,y)z_2(x,y) = -\left( \tfrac{x}{a} \right)^2 -  \left( \tfrac{y}{b} \right)^2.
\end{align*}
These surfaces describe the asymptotic behavior near conical cusps. Unfortunately the problem in (\ref{eq:toyconical}--\ref{eq:toyconicalend}) is too simple for a computational toy problem since $s_1, s_2$ are low degree polynomials and the reconstructions using companion matrix methods are thus almost immediately exact to their full precision. Thus, to make the problem more interesting, we consider a non-algebraic modification
\begin{align}\label{eq:toyconicalsinh}
    \bar{z}_1(x,y) &= \sinh\left(\tfrac{\sqrt{a^2 y^2+b^2 x^2}}{a b}\right),\\
    \bar{z}_2(x,y) &= \sinh\left(-\tfrac{\sqrt{a^2 y^2+b^2 x^2}}{a b}\right),\label{eq:toyconicalsinhend}
\end{align}
which also has a conical cusp at $(x,y) = (0,0)$ with smooth ESPs:
\begin{align*}
    \bar{s}_1(\bar{\mathbf{z}}(x,y)) &= \bar{z}_1(x,y) + \bar{z}_2(x,y) = 0,\\
    \bar{s}_2(\bar{\mathbf{z}}(x,y)) &= \bar{z}_1(x,y)\bar{z}_2(x,y) = -\sinh^2\left(\tfrac{\sqrt{a^2 y^2+b^2 x^2}}{a b}\right).
\end{align*}
We plot the surfaces $z_1(x,y), z_2(x,y)$ as well as the smooth ESP $\bar{s}_2(\mathbf{z}(x,y))$ with $(a,b) = (\frac{4}{3}, \frac{12}{5})$ in Fig. \ref{fig:sinhconicaltoyexampledemo} and the error of the Chebyshev colleague approach for both problems in Fig. \ref{fig:2dsinusoidaltoyfigreconstruction}(b), showing that the method performs well for conical cusps. In Section \ref{sec:graphene} we explore the electronic band structure of graphene which is an example of an application featuring multiple conical cusps.
\subsubsection{Non-smooth ESPs}
If all to-be-reconstructed surfaces are smooth then it is guaranteed that the ESPs as well as the Chebyshev equivalents will be smooth. Furthermore, we have seen in the previous section that for certain cases of interest, e.g. conical intersections, it is possible that the ESPs are smooth even if the underlying surfaces are not. In this section we consider a special case relevant for applications featuring non-smooth ESPs.

For simplicity we will consider a 1D example: We begin with a stack of $m$ surfaces, ordered by their lowest attained value such that every surface intersects the next at least once in the domain of interest. We then consider the problem of reconstructing $m-1$ underlying surfaces from value ordered data such that at least one non-differentiable cusp remains in the highest value entry of our multi-surface data. Specifically, we will be considering the same problem from Section \ref{sec:numericaltoy1} and Figure \ref{fig:firsttoyexample} but adding a fourth, not-to-be-reconstructed surface
\begin{align}\label{eq:additionalforstack}
    z_4(x) = \tfrac{1}{3}+\cos(\tfrac{2x}{3}).
\end{align}
This scenario is significant for applications since e.g. the hierarchy of excited states involves in principle infinitely many stacked energy surfaces where one must choose an application-driven but ultimately arbitrary cutoff point.

\begin{figure}\centering
     \subfloat[]
    {{ \centering \includegraphics[height=6.7cm]{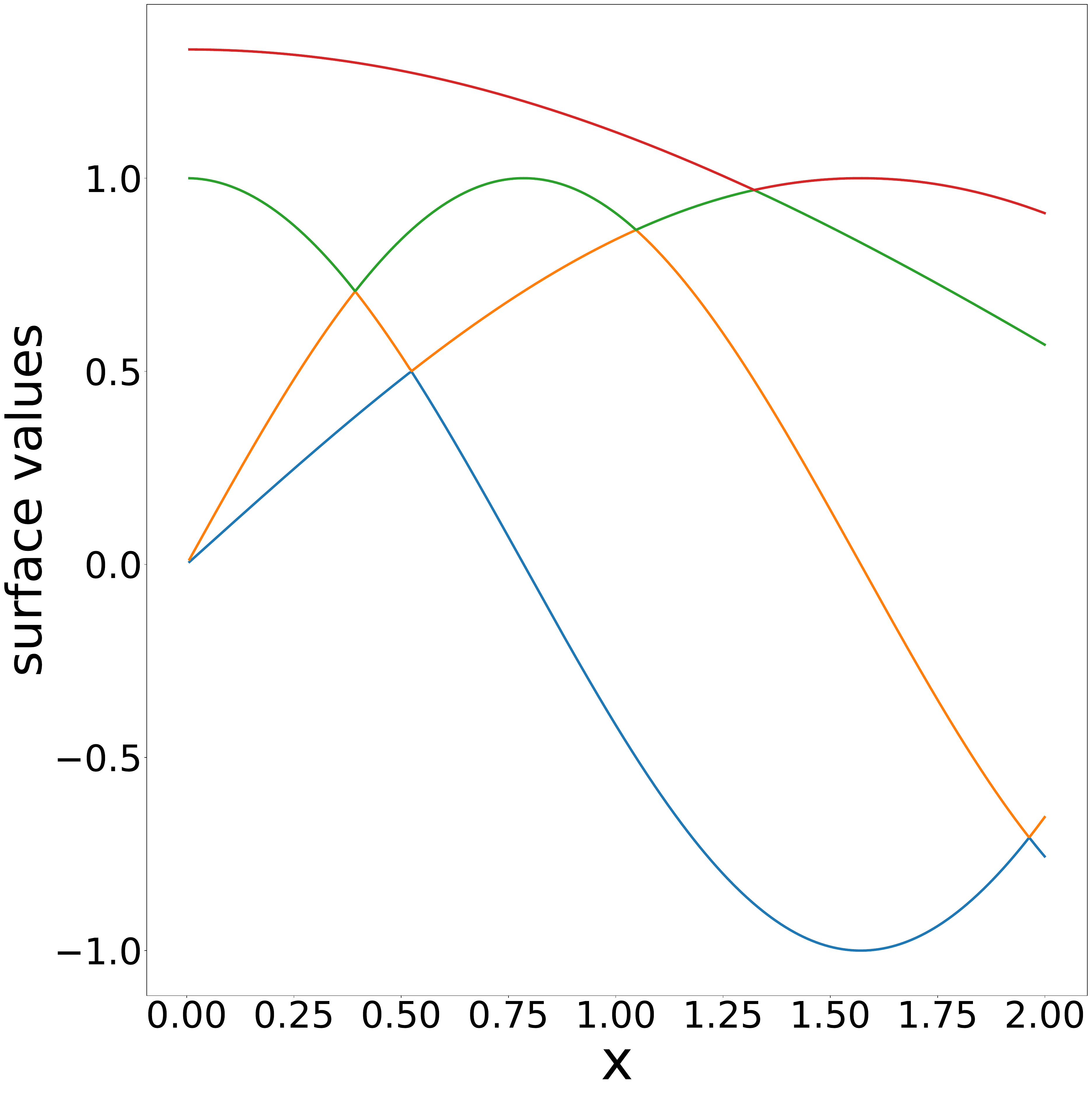} }}
     \subfloat[]
    {{ \centering \includegraphics[height=6.7cm]{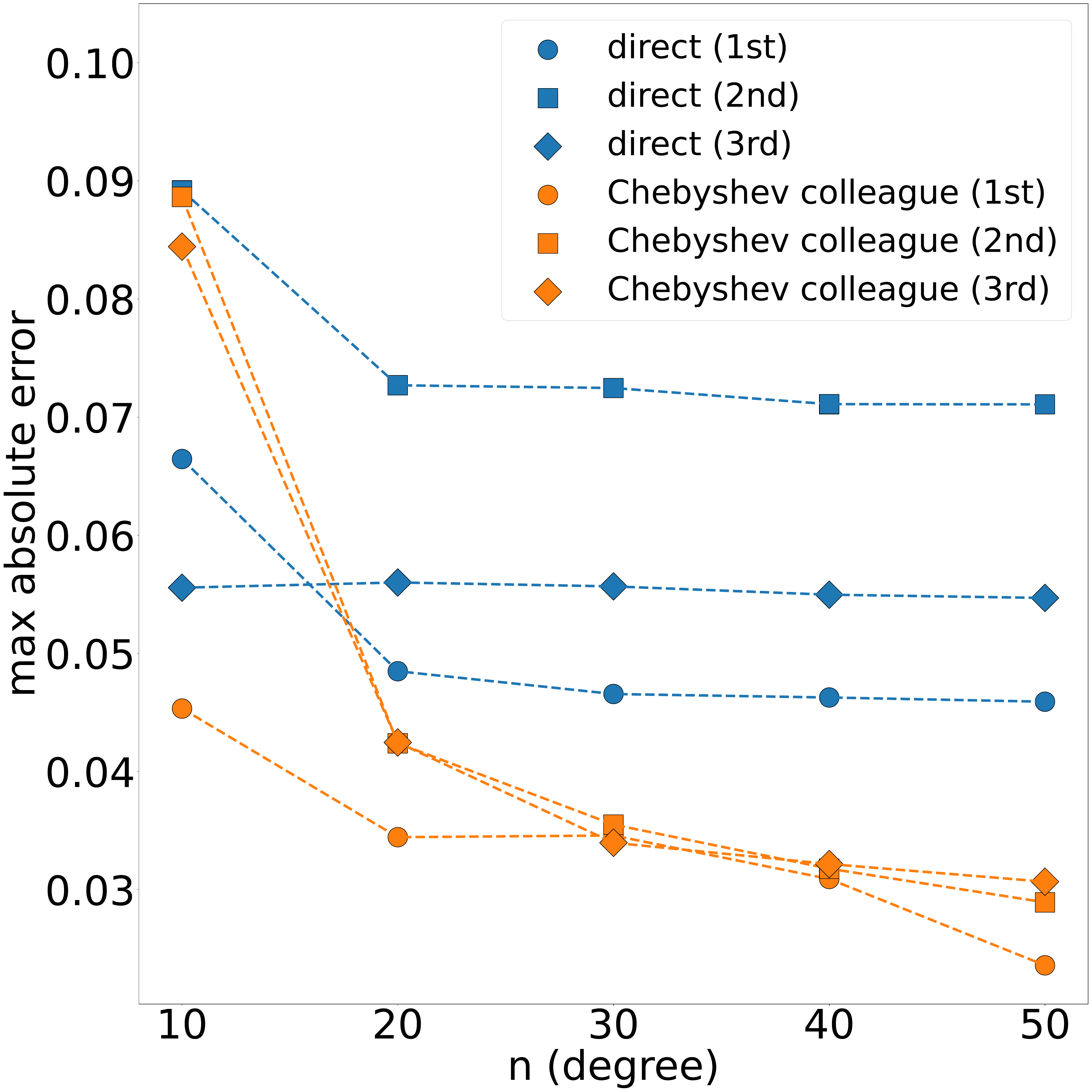} }}
    \caption{(a) Intersecting, smooth sinusoidal curves in (\ref{eq:toy1start}--\ref{eq:toy1end}) and \eqref{eq:additionalforstack} pointwise sorted by values, cf. Figure \ref{fig:firsttoyexample}. Reconstructing only the lowest three layers leaves a cusp in the third layer near $x\approx 1.32$. (b) Linear scale max. abs. error plot for reconstructing the lowest three value-sorted surfaces from (a).}
    \label{fig:stackedsurfacecase}
\end{figure}

We plot max. abs. errors for reconstructing the first three value-sorted surfaces described in (\ref{eq:toy1start}--\ref{eq:toy1end}) and \eqref{eq:additionalforstack} using direct interpolation and Chebyshev colleague matrix methods in Figure \ref{fig:stackedsurfacecase}(b). One observes that while the colleague matrix approach outperforms the direct method, the differences are less pronounced than in the case of smooth invariants. The observed behavior is expected since including a cusp in the invariants causes similar problems for the companion matrix approaches as it does for the direct method. The fact we still observe better errors is likewise expected since the number of cusps in the invariants is less than or equal to the number of unmatched / leftover cusps in the top surface times the number of surfaces to reconstruct, which is always less than or equal to the number of cusps in the original data. Informally, fewer cusps in the data lead to an expectation of better polynomial approximations, albeit only by constant factors. The companion matrix methods perform best in stacked surface scenarios if the cutoff surface is chosen such that the fewest possible number of cusps are left unmatched in the top layer. It is an interesting question whether there are practical ways to artificially match the leftover cusp in the top layer to restore smooth invariants.

\subsection{Example applications}
\subsubsection{Conical cusps in the electronic band structure of graphene}\label{sec:graphene}
Graphene consists of carbon atoms arranged in a hexagonal honeycomb lattice. For a single layer of graphene a tight-binding model (considering only nearest-neighbor interactions) can be used to obtain a Hamiltonian whose conduction and valance bands can be explicitly computed \cite{liu2018graphene,semenoff1984condensed,wallace1947band}:
\begin{align}\label{eq:grapheneconductandvalence}
E(\mathbf{k})=\pm\,\gamma_0\sqrt{1+4\cos^2\left({\tfrac{a k_x}{2}}\right)+4\cos\left({\tfrac{a k_x}{2}}\right) \cos\left({\tfrac{\sqrt{3}a k_y}{2}}\right)},
\end{align}
where $\mathbf{k} = [k_x, k_y]$ denotes the electron wave vector and $\gamma_0 \approx 2.8 \text{eV}$, $a \approx 2.46 \text{\AA}$ are physical constants. For thorough discussions of the physics and chemistry of graphene we refer to \cite{neto2009electronic,liu2018graphene} and the references therein.

\begin{figure}\centering
      \subfloat[]
    {{ \centering \includegraphics[width=7.3cm]{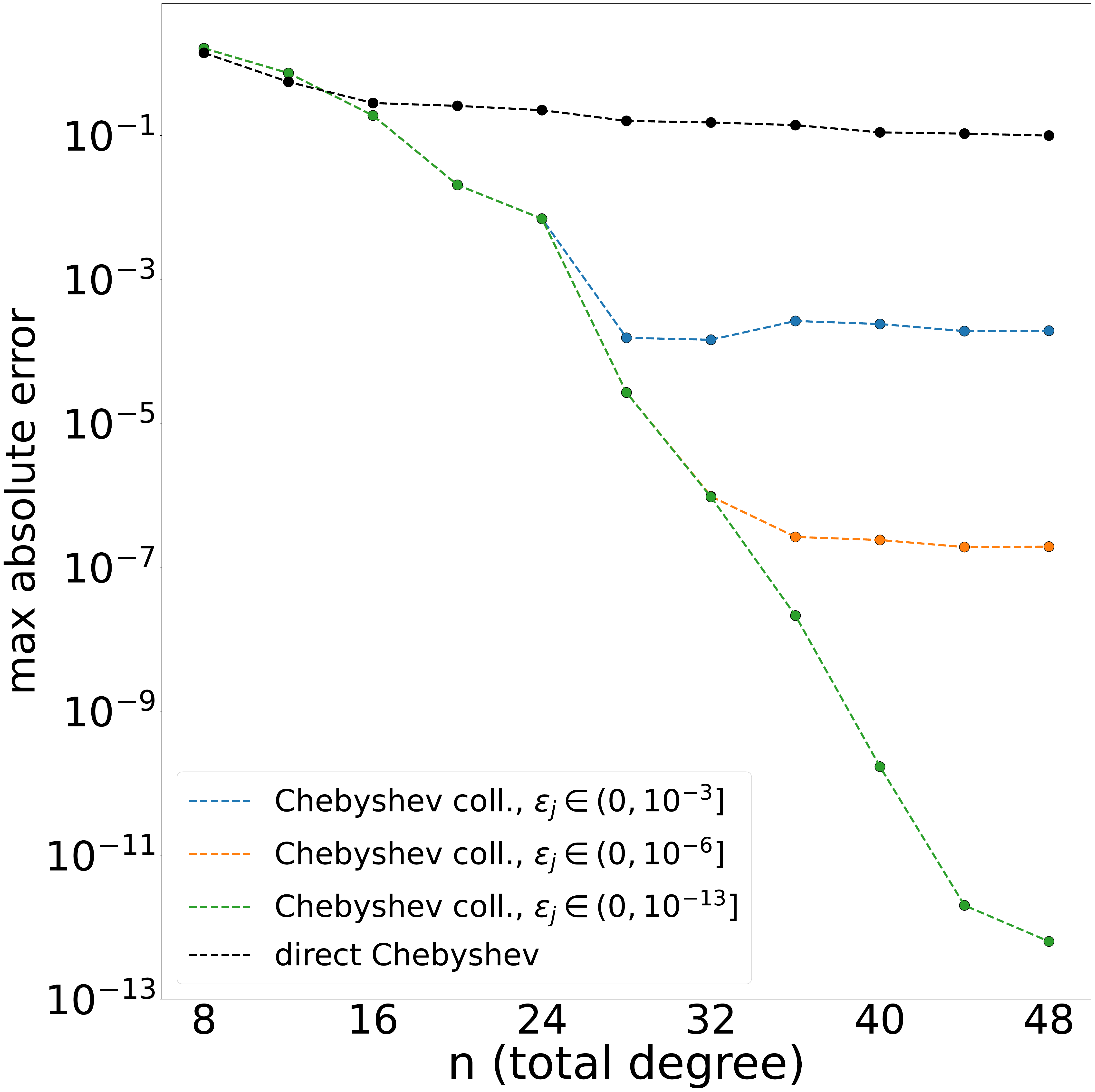}}}
     \subfloat[]
    {{ \centering \includegraphics[width=7.3cm]{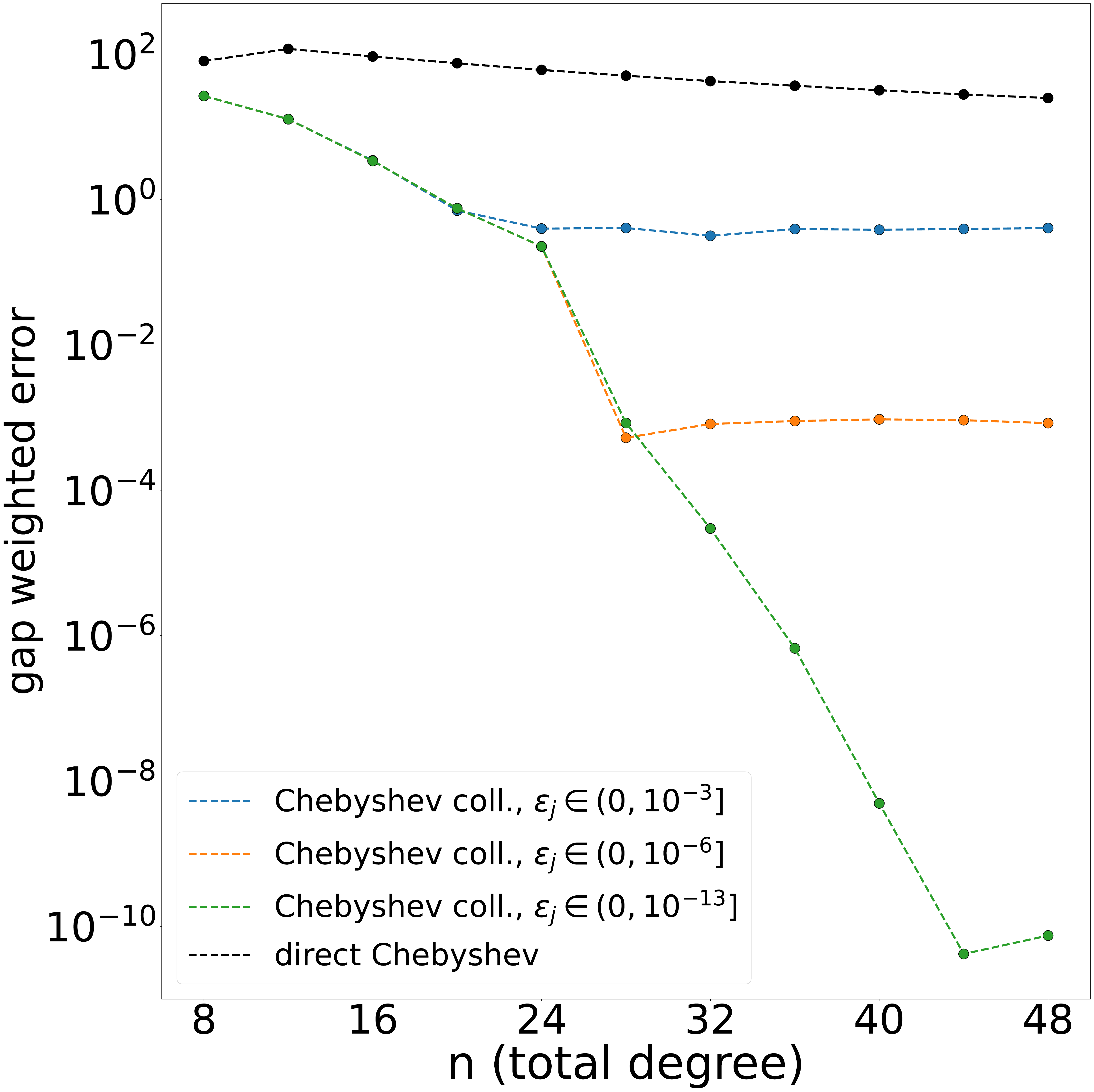}}}
    \caption{Max. abs. error in (a) and gap weighted error in (b) on a $150 \times 150$ grid using the Chebyshev colleague matrix and direct interpolation methods with varying data perturbation $\epsilon_j$ for the Graphene example described in \eqref{eq:grapheneconductandvalence}.}
    \label{fig:graphenemain}
\end{figure}

\begin{figure}\centering
     \subfloat[]
    {{ \centering \includegraphics[width=7.3cm]{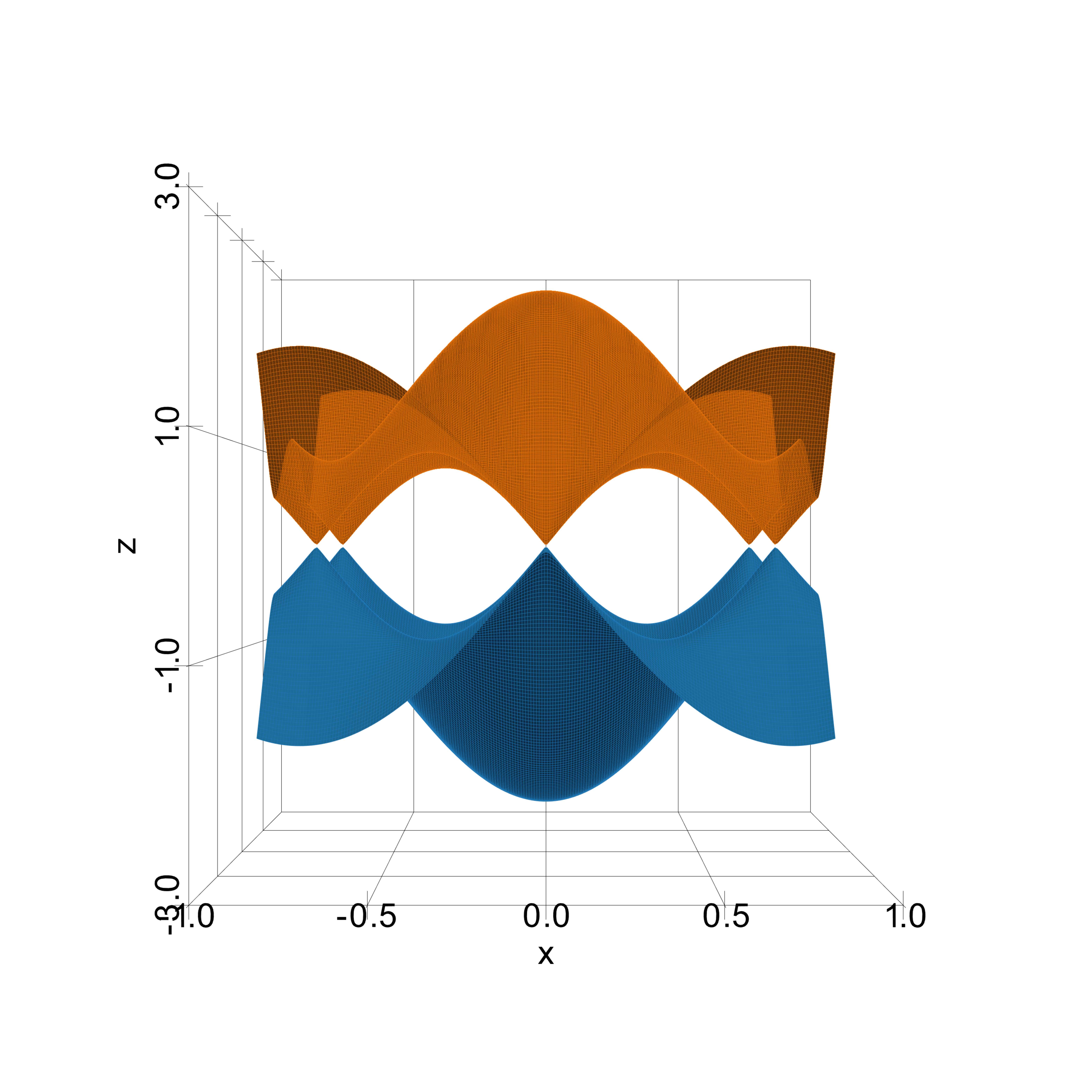} }}
     \subfloat[]
    {{ \centering \includegraphics[width=7.3cm]{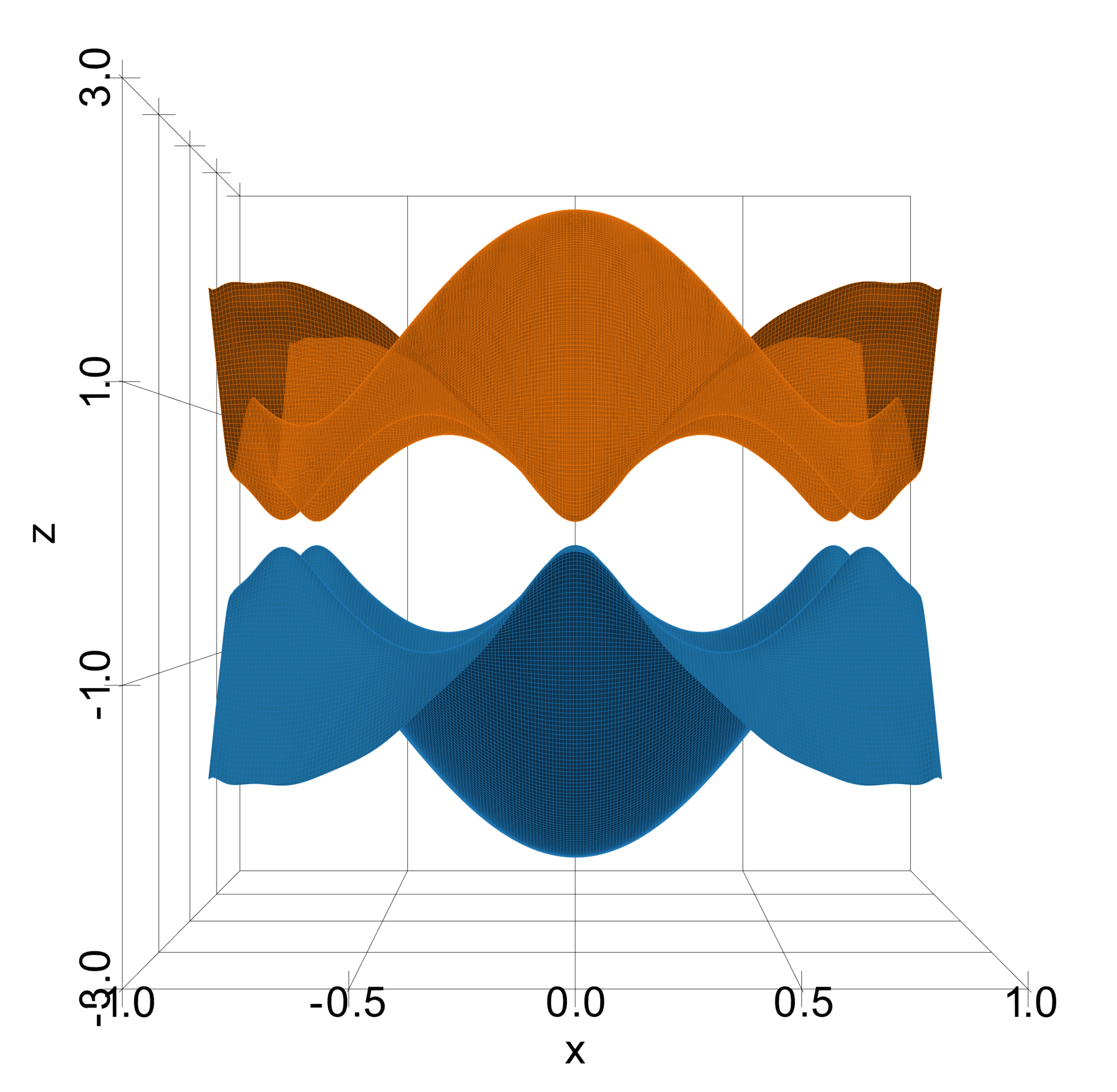} }}
    \caption{(a) shows scaled, reconstructed Graphene conduction and valence bands using a Chebyshev colleague method (tot. degree $n=40$) with the underlying data randomly perturbed by $\epsilon_j \in (0,10^{-3}]$, (b) shows analogous reconstruction using direct approximation with Chebyshev polynomials (tot. degree $n=40$) and exact data.}
    \label{fig:visualcompgraphene}
\end{figure}

As seen in Figure \ref{fig:visualcompgraphene}(a), graphene's conduction and valence bands touch at the so-called Dirac points, making it a zero-gap semiconductor. In the language of the previous sections this means that the conduction-valence multi-surface of graphene has multiple conical intersections. The electronic band structure of graphene thus constitutes a more advanced and application-oriented version of results discussed in Section \ref{sec:numericaltoy2cusp}. In Figure \ref{fig:graphenemain} we plot the maximum absolute error for reconstructing the bands with various levels of perturbations prior to polynomial interpolation of the invariants.

In Figure  \ref{fig:visualcompgraphene} we provide a visual comparison of the avoided crossings obtained using a direct interpolation of the surfaces versus those obtained via the Chebyshev colleague method, showing substantially better resolved Dirac cones. Figure \ref{fig:visualcompgraphene} also provides a clear demonstration of another key strength of this approach: Unlike a direct interpolation attempt, the presence of cusps does not pollute the rest of the multi-surface.
\subsubsection{Energy surfaces of a sulfur dioxide ($SO_2$) molecule}\label{sec:so2}
\begin{figure}
 \centering \includegraphics[width=7.3cm]{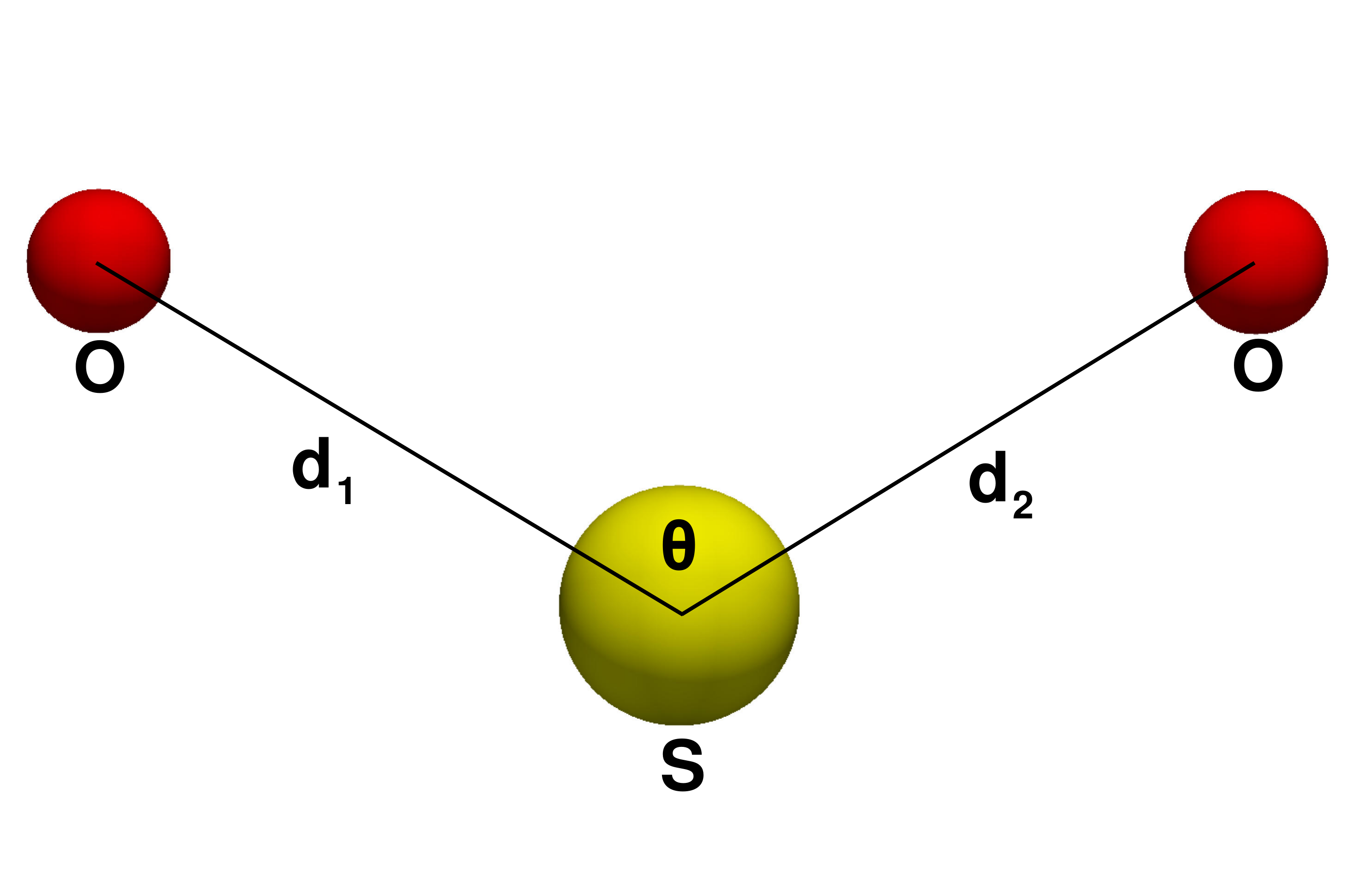}
    \caption{Schematic of $SO_2$ molecule described by coordinate triplet $(\theta, d_1, d_2)$.}
    \label{fig:so2}
\end{figure}

In this section we use our method to reconstruct the energy surfaces for sulfur dioxide ($S O_2$), a molecule consisting of a sulfur and two oxygen atoms as pictured in Figure \ref{fig:so2}, from global interpolants. We pre-compute (unordered) energies associated with different relative positions for different angles $\theta$ and distances $d_1$ and $d_2$ between the sulfur atom and respective oxygen atoms using a linear vibronic coupling (LVC) model of SO$_{2}$ \cite{westermayr_combining_2020,plasser2019highly, kouppel1984multimode}. 

For larger molecules it would be natural to machine learn the ESPs from limited pre-computed samples using frameworks such as the atomic cluster expansion (ACE) \cite{DrautzACE,ACECompleteness} or MACE \cite{Batatia2022mace}. The example of $S O_2$ is sufficiently small, however, that we can treat the problem directly using Chebyshev polynomials since as seen in Figure \ref{fig:so2} the state is up to global symmetries fully determined from coordinate triplets such as $(\theta, d_1, d_2)$ or $(\theta, d_1-d_2, d_1+d_2)$. Since a conical intersection is known to occur in the potential energy surfaces for $d_1-d_2 \approx 0$, cf. \cite{wilkinson2014excited}, we will use the latter set of internal coordinates.

\begin{figure}\centering
\subfloat[]{{
 \centering \includegraphics[width=7.3cm]{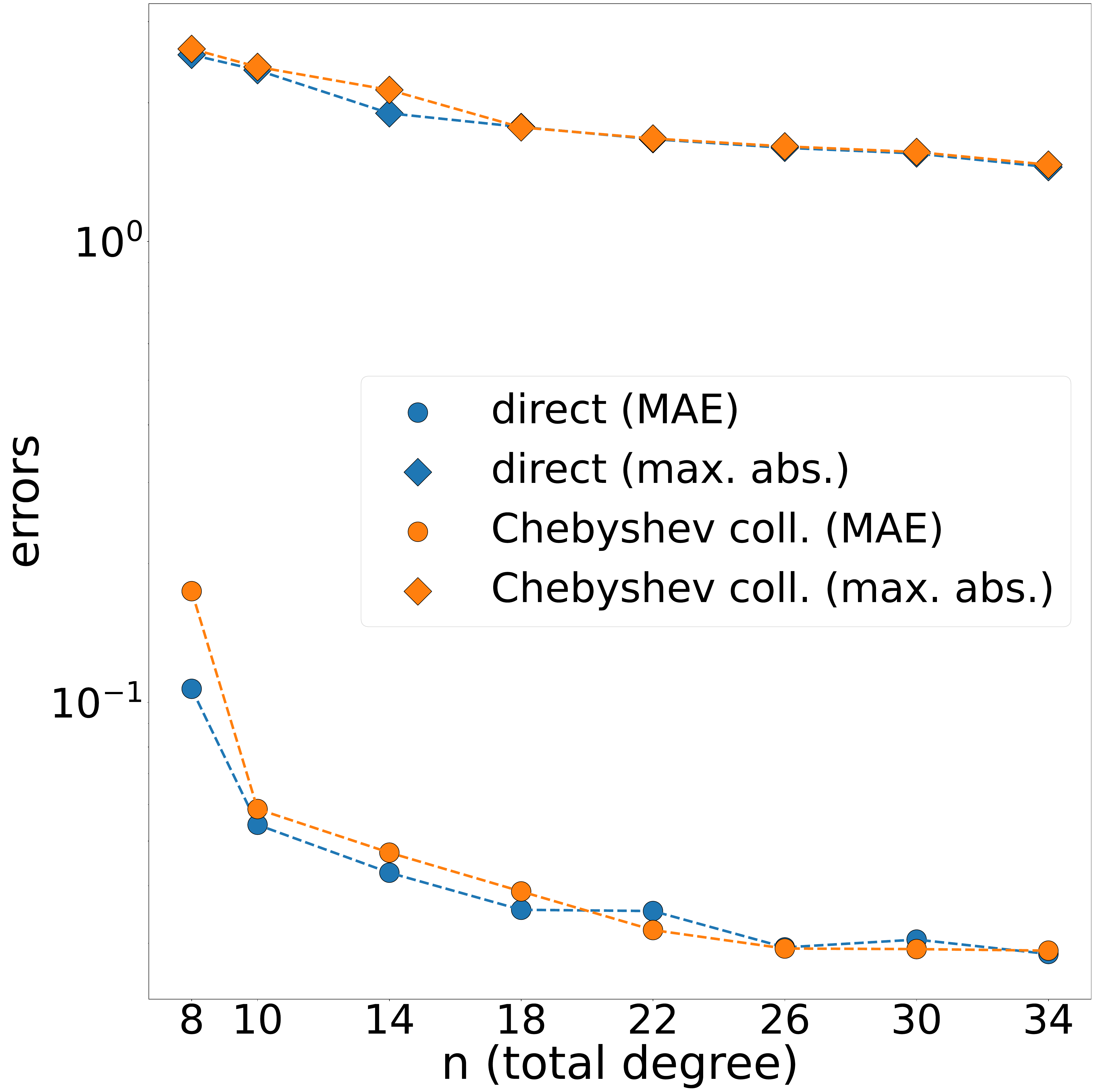}}}
 \subfloat[]{{
 \centering \includegraphics[width=7.3cm]{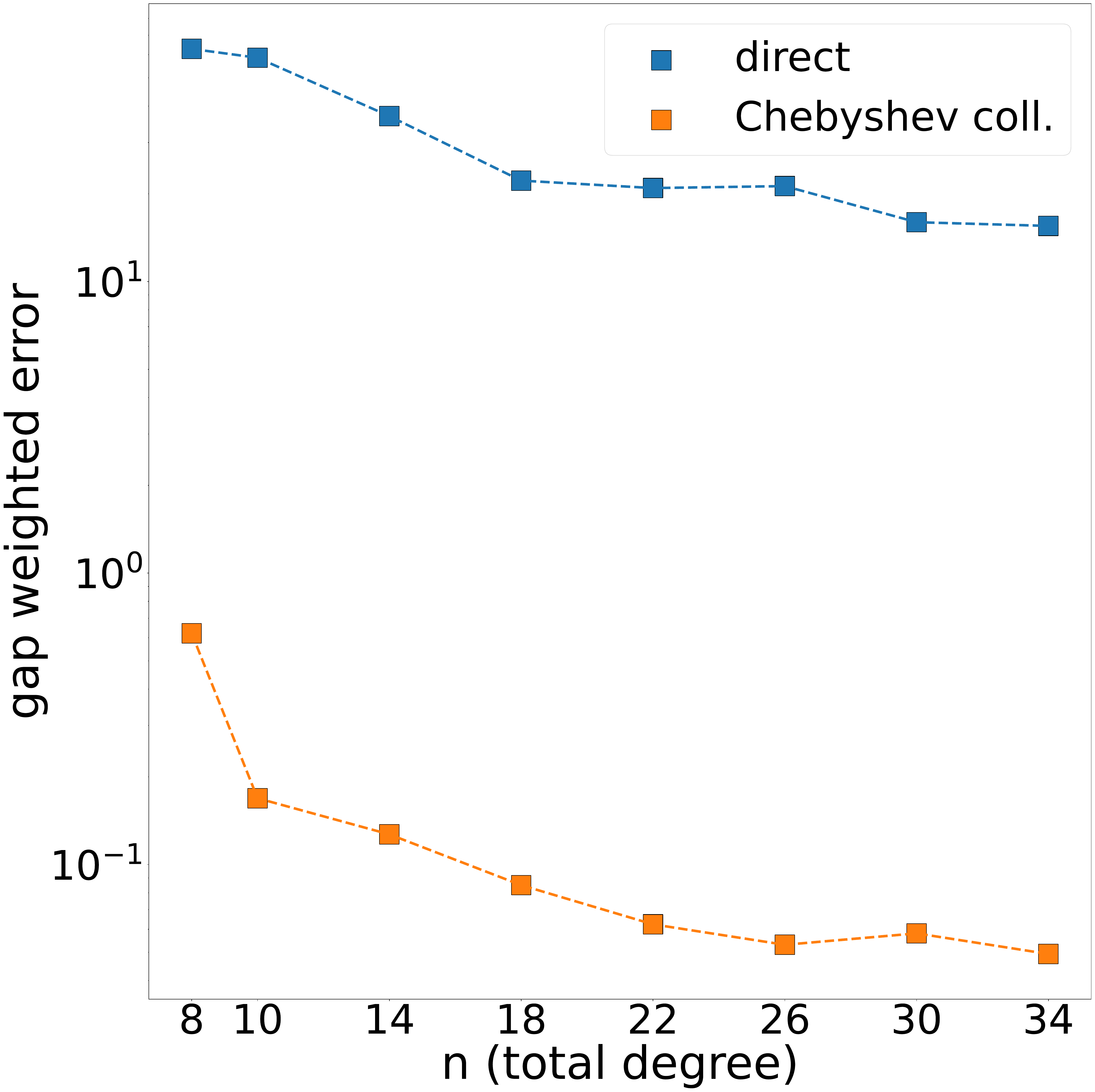}}}
    \caption{(a) Semi-logarithmic error plot showing obtained MAE and max. abs. errors of the direct and Chebyshev colleague matrix methods associated with different configurations of $SO_2$, cf. Fig. \ref{fig:so2cuspslices}. (b) is an analogous plot showing gap weighted errors.}
    \label{fig:so2errors}
\end{figure}

\begin{figure}\centering
\subfloat[cusp in $SO_2$ data projection]
{{\centering
    \includegraphics[width=5.0cm]{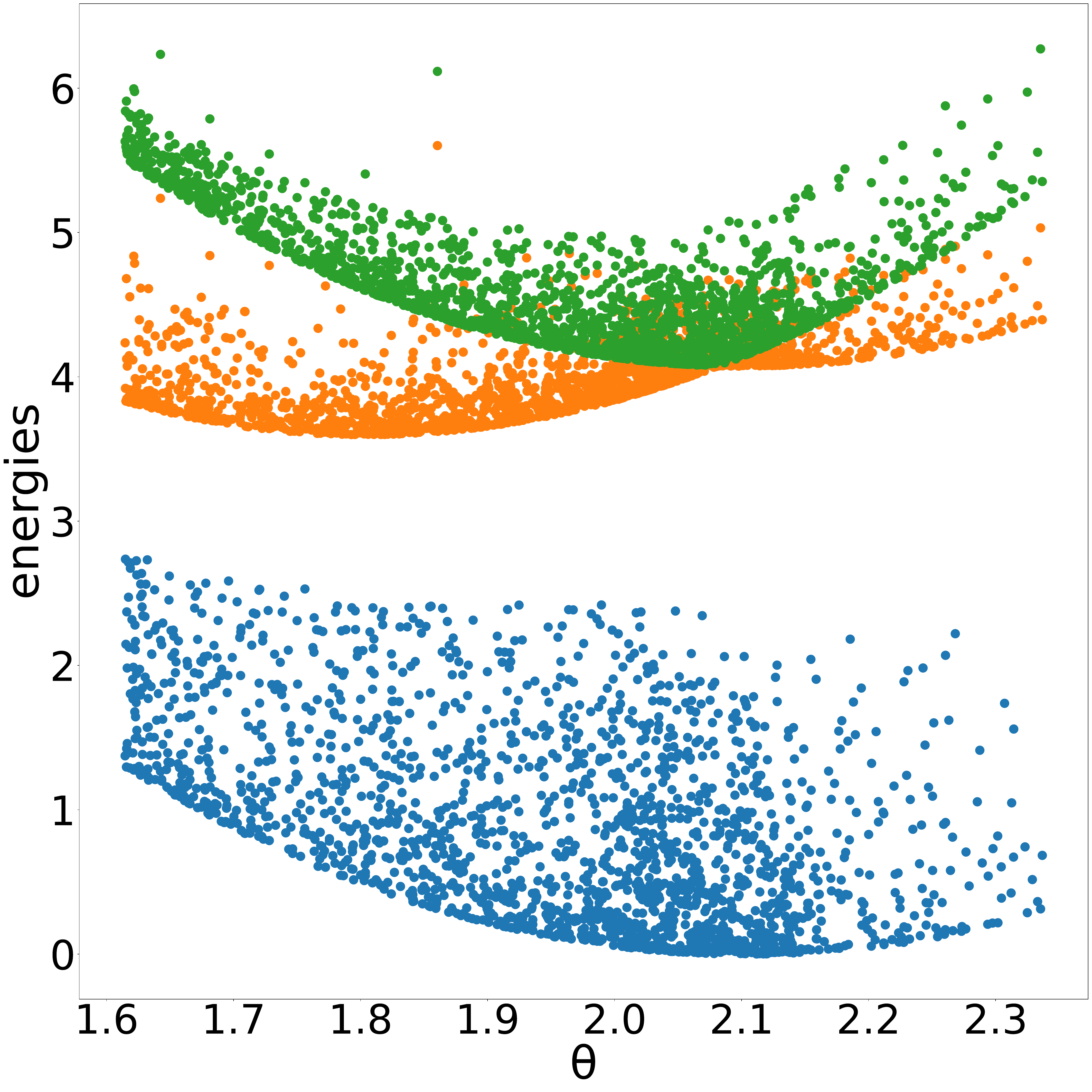}}}
    \subfloat[Chebyshev colleague ($n=28$)]
{{\centering
    \includegraphics[width=5.0cm]{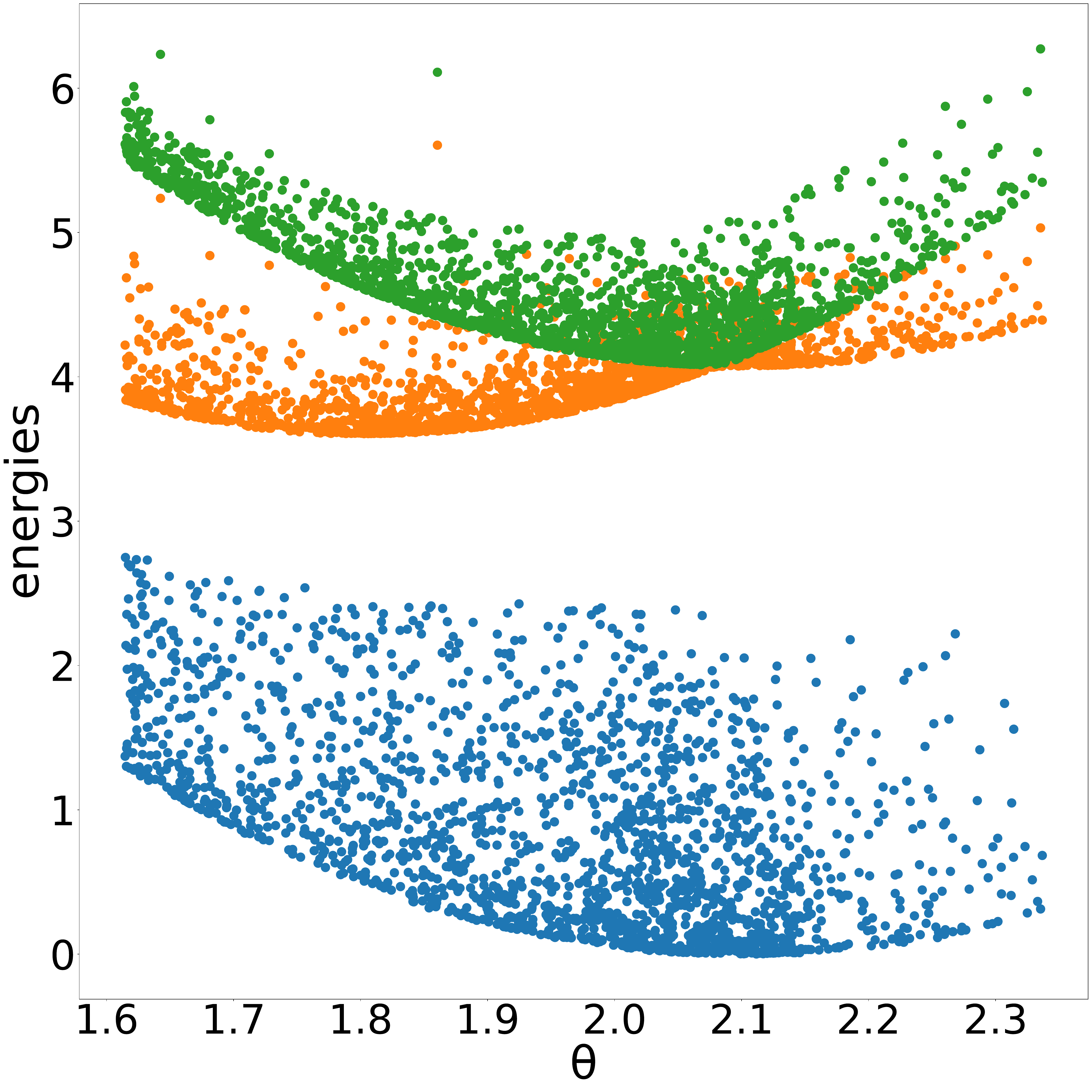}}}
    \subfloat[direct method ($n=28$)]
{{\centering
    \includegraphics[width=5.0cm]{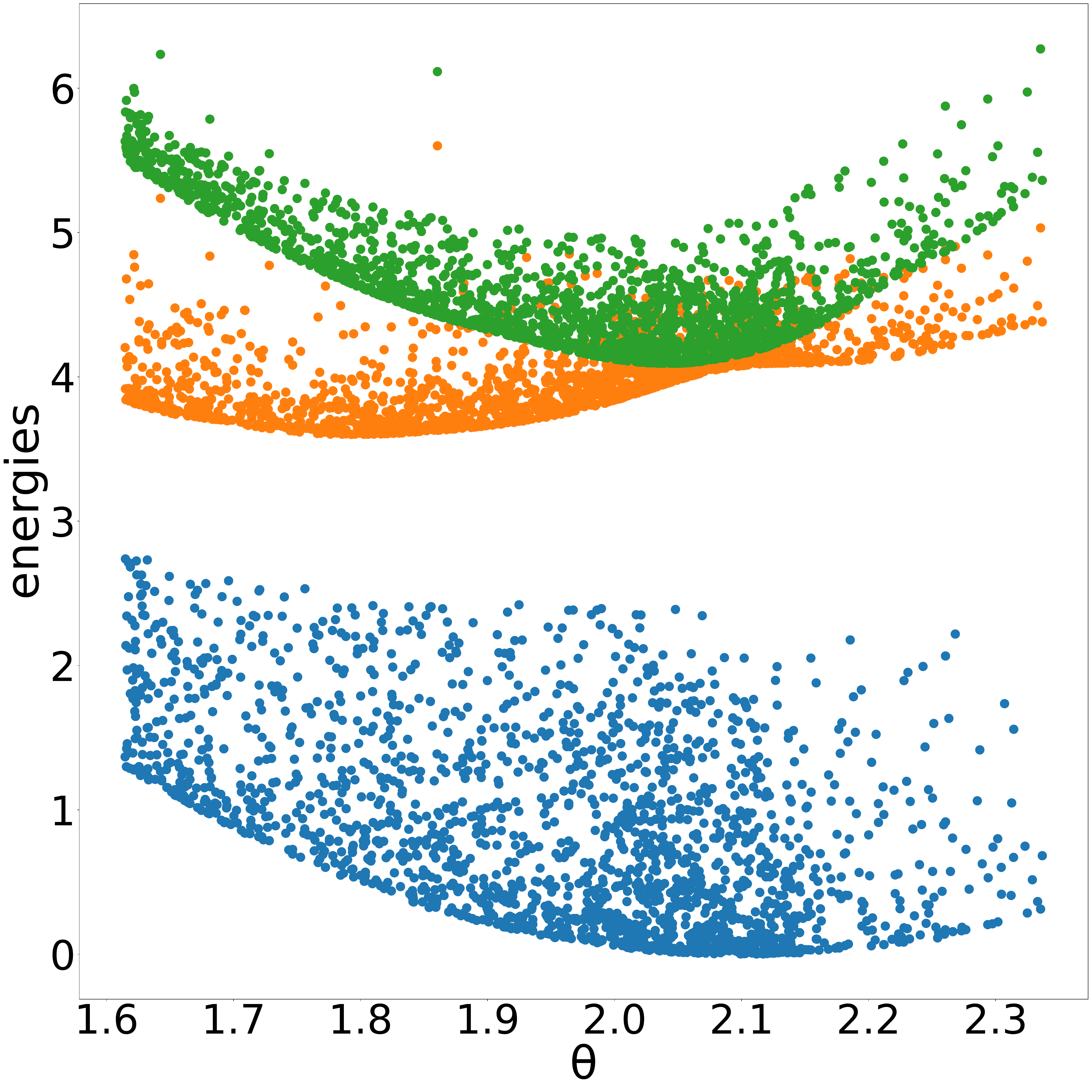}}}\\
    \subfloat[zoom in of (a)]
    {{\centering
    \includegraphics[width=5.0cm]{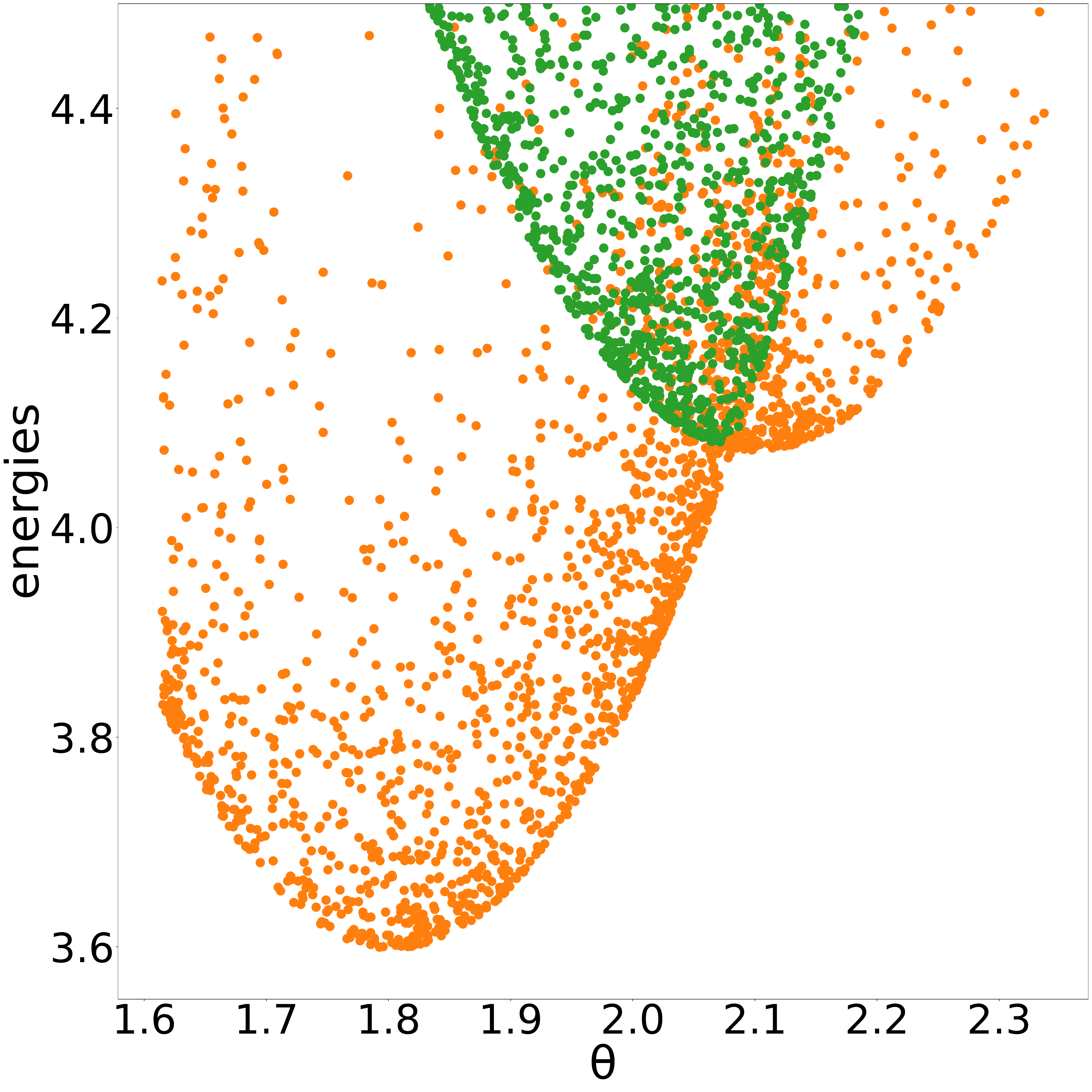}}}
    \subfloat[zoom-in of (b)]
{{\centering
    \includegraphics[width=5.0cm]{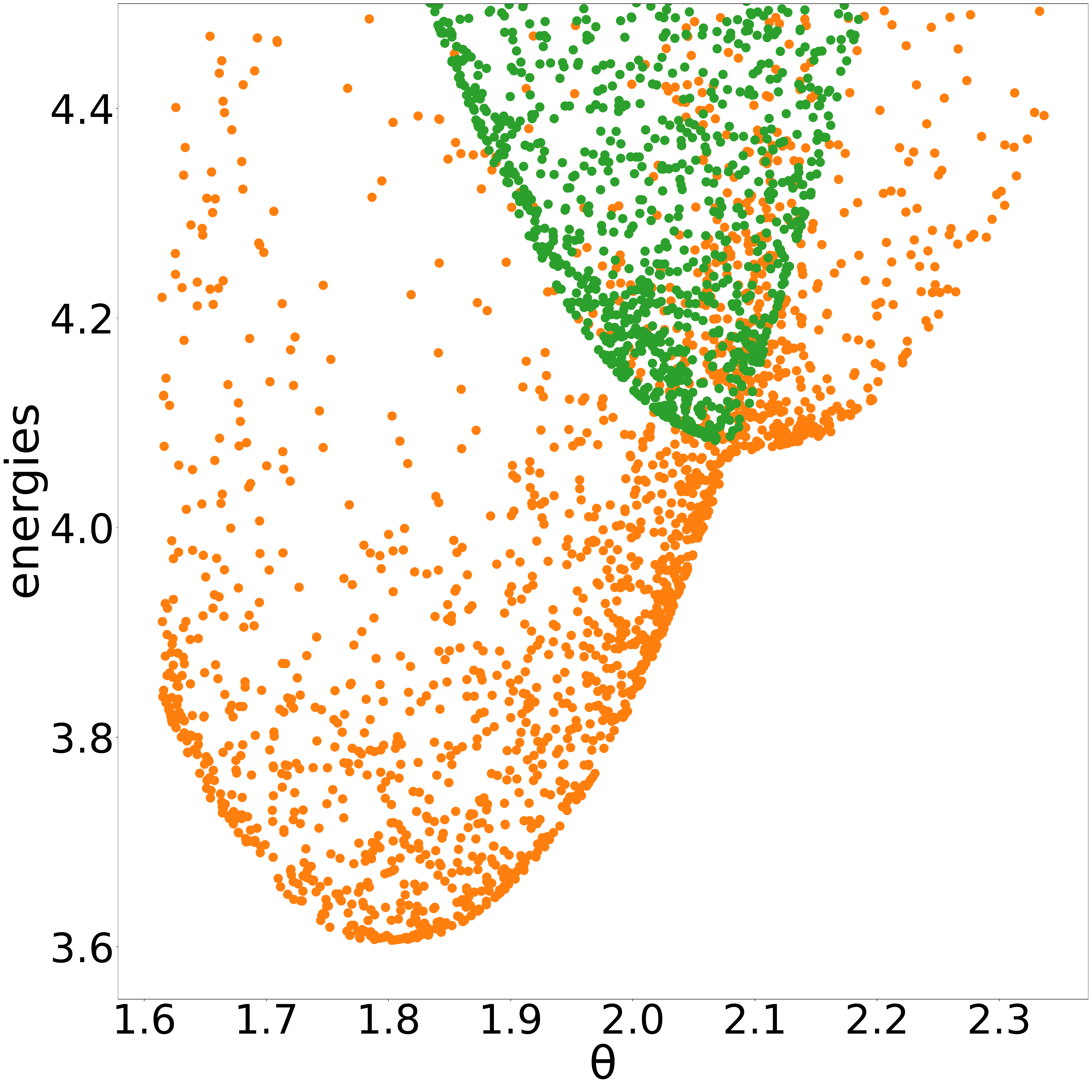}}}
    \subfloat[zoom-in of (c)]
{{\centering
    \includegraphics[width=5.0cm]{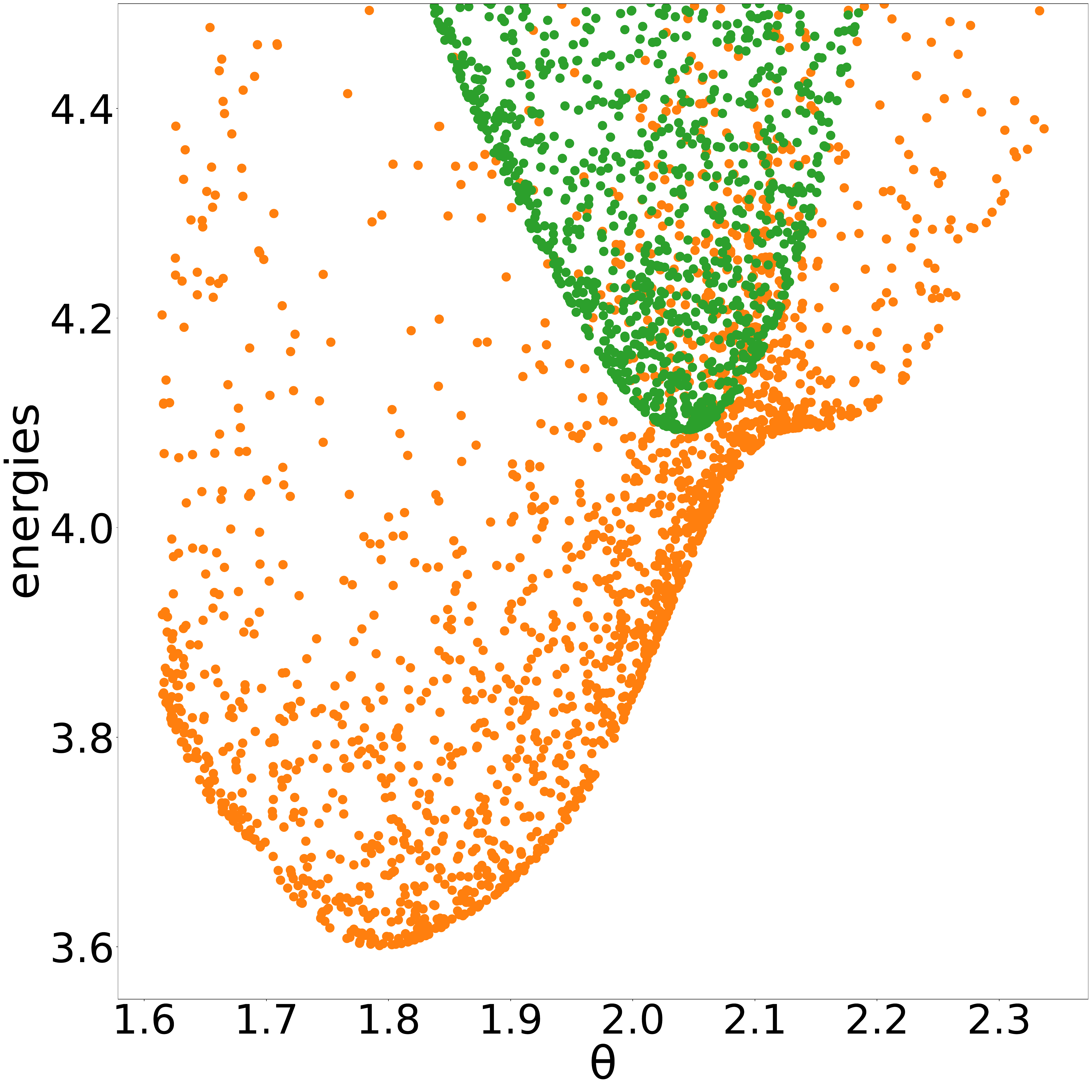}}}
    \caption{(a) shows cusp presence for a three-dimensional $SO_2$ dataset. The presented scatter plots are 1D $\theta$-projections of 2D-slices with $d_1-d_2 \approx 0$ aiming to visualize the presence of a cusp. (b) and (c) show reconstructions using Chebyshev colleague and direct methods at total degree $n=28$. (d-f) show zoomed-in segments of (a-c).}
    \label{fig:so2cuspslices}
\end{figure}

As the energy surfaces for this setup are three-dimensional, a direct visualization is effectively impossible. We thus instead plot projected slices of our data in Figure \ref{fig:so2cuspslices} showing the presence of a cusp. Since the underlying data is obtained from an LVC model, giving approximations to the actual electronic structure of the system, the results of this section are to be interpreted as relying on highly perturbed data of a non-negligible but unspecified amount similar to the scenario one would encounter in more complex computational chemistry applications. Errors due to outliers in the data dominate the max. abs. errors as well as overall MAE when reconstructing the surfaces, see Figure \ref{fig:so2errors}(a), and thus as discussed at the beginning of Section \ref{sec:experiments} mask the improvement near the cusps plainly visible in Figure \ref{fig:so2cuspslices}. The gap weighted error defined in Eq. \eqref{eq:gapweightederror} which we plot in Figure \ref{fig:so2errors}(b) provides a natural quantitative measure to distinguish between the cusp resolution quality of direct reconstructions and the Chebyshev colleague matrix approach.

\section{Discussion}
We described a Schmeisser companion and Chebyshev colleague method to reconstruct a multi-surface from unordered or value-sorted data and explored their numerical properties. The methods have attractive approximation properties in low noise context and retain some of these properties into the more application-relevant noisy regimes. The primary expected value of these methods is improved resolution of multi-surface cusps with globally smooth interpolations or alternatively equivalent quality resolution of cusps at lower polynomial degree.

Throughout the paper we have mentioned opportunities to add pre- or post-processing steps in order to further improve the practical usefulness of these methods. An alternative approach could be the use of different smooth invariants altogether. A particularly interesting future direction would be the combination of such companion matrix based methods with machine learning models for molecular properties such as the atomic cluster expansion (ACE) \cite{DrautzACE,ACECompleteness} and related MACE \cite{Batatia2022mace} frameworks which rely on polynomial interpolations and are thus a natural fit.

\section{Acknowledgements}
This work was supported by NSERC Discovery Grant GR019381 and NFRF Exploration Grant GR022937. 
TSG was also supported by a PIMS-Simons postdoctoral fellowship, jointly funded by the Pacific Institute for the Mathematical Sciences (PIMS) and the Simons Foundation. This work was also supported by the Deutsche Forschungsgemeinschaft (DFG) – Project-ID 443871192 - GRK 2721: "Hydrogen Isotopes 1,2,3H." The authors would like to acknowledge the ZIH TU Dresden and the URZ Leipzig University for providing computational resources.

The authors acknowledge helpful conversations with Sebastian Mai and Cheuk Hin Ho on implementations, applications and motivation of the presented methods as well as Isaac Holt for feedback on the manuscript.

\newpage
\printbibliography[heading=bibintoc]

@article{hladik2017eigenvalues,
  title={Eigenvalues of symmetric tridiagonal interval matrices revisited},
  author={Hlad{\'\i}k, M.},
  journal={arXiv preprint arXiv:1704.03670},
  year={2017}
}

@Article{golub1962bounds,
  author    = {Golub, G. H},
  journal   = {Math. Comput.},
  title     = {Bounds for eigenvalues of tridiagonal symmetric matrices computed by the LR method},
  year      = {1962},
  number    = {80},
  pages     = {438--445},
  volume    = {16},
  fjournal  = {Mathematics of Computation},
  publisher = {JSTOR},
}

@Article{wilkinson1958calculation,
  author    = {Wilkinson, J. H.},
  journal   = {Comput. J.},
  title     = {The calculation of eigenvectors by the method of Lanczos},
  year      = {1958},
  number    = {3},
  pages     = {148--152},
  volume    = {1},
  fjournal  = {The Computer Journal},
  publisher = {Oxford University Press},
}

@book{parlett1998symmetric,
  title={The symmetric eigenvalue problem},
  author={Parlett, B. N.},
  year={1998},
  publisher={SIAM}
}

@Article{hladik2011characterizing,
  author    = {Hlad{\'\i}k, M. and Daney, D. and Tsigaridas, E.},
  journal   = {Comput. Math. Appl.},
  title     = {Characterizing and approximating eigenvalue sets of symmetric interval matrices},
  year      = {2011},
  number    = {8},
  pages     = {3152--3163},
  volume    = {62},
  fjournal  = {Computers & Mathematics with Applications},
  publisher = {Elsevier},
}

@Article{hladik2010bounds,
  author    = {Hlad{\'\i}k, M. and Daney, D. and Tsigaridas, E.},
  journal   = {SIAM J. Matrix Anal. Appl.},
  title     = {Bounds on real eigenvalues and singular values of interval matrices},
  year      = {2010},
  number    = {4},
  pages     = {2116--2129},
  volume    = {31},
  fjournal  = {SIAM Journal on Matrix Analysis and Applications},
  publisher = {SIAM},
}

@article{rohn2005handbook,
  title={A handbook of results on interval linear problems},
  author={Rohn, Jir{\i}},
  year={2005},
journal = {Technical report No. V-1163, Institute of Computer Science, Academy of Sciences of the Czech Republic}
}

@Article{schmeisser1993real,
  author    = {Schmeisser, G.},
  journal   = {Linear Algebra Appl.},
  title     = {A real symmetric tridiagonal matrix with a given characteristic polynomial},
  year      = {1993},
  pages     = {11--18},
  volume    = {193},
  publisher = {Elsevier},
}

@Article{fiedler1990expressing,
  author    = {Fiedler, M.},
  journal   = {Linear Algebra Appl.},
  title     = {Expressing a polynomial as the characteristic polynomial of a symmetric matrix},
  year      = {1990},
  pages     = {265--270},
  volume    = {141},
  fjournal  = {Linear Algebra and its Applications},
  publisher = {Elsevier},
}

@Article{fiedler2003note,
  author    = {Fiedler, M.},
  journal   = {Linear Algebra Appl.},
  title     = {A note on companion matrices},
  year      = {2003},
  pages     = {325--331},
  volume    = {372},
  fjournal  = {Linear Algebra and its Applications},
  publisher = {Elsevier},
}

@book{smith1985numerical,
  title={Numerical solution of partial differential equations: finite difference methods},
  author={Smith, G. D.},
  year={1985},
  publisher={Oxford University Press}
}

@Article{noschese2013tridiagonal,
  author    = {Noschese, S. and Pasquini, L. and Reichel, L.},
  journal   = {Numer. Linear. Algebr. Appl.},
  title     = {Tridiagonal Toeplitz matrices: properties and novel applications},
  year      = {2013},
  number    = {2},
  pages     = {302--326},
  volume    = {20},
  publisher = {Wiley Online Library},
}

@book{egge2019introduction,
  title={An introduction to symmetric functions and their combinatorics},
  author={Egge, E. S.},
  volume={91},
  year={2019},
  publisher={American Mathematical Soc.}
}

@book{Stanley_Fomin_1999, place={Cambridge}, series={Cambridge Studies in Advanced Mathematics}, title={Enumerative Combinatorics}, publisher={Cambridge University Press}, author={Stanley, R. P. and Fomin, S.}, year={1999}, collection={Cambridge Studies in Advanced Mathematics}}

@book{loehr2017combinatorics,
  title={Combinatorics},
  author={Loehr, Nicholas},
  year={2017},
edition = {Second Edition},
  publisher={Chapman and Hall/CRC}
}

@book{vietae1646,
  author    = {Vieta, F.},
  title     = {Opera Mathematica},
  editor    = {Van Schooten, Frans},
  publisher = {Elzévir},
  address   = {Leyde},
  year      = {1646},
  note      = {Hildesheim-New-York: Georg Olms Verlag, 1970 (in Latin)}
}

@Book{girard1884,
  author  = {Girard, A.},
  editor  = {Bierens de Haan, D.},
  title   = {Invention Nouvelle En L'algebre (Reprint)},
  year    = {1884},
  address = {Leiden},
  note    = {(in French)},
}

@Article{funkhouser1930short,
  author    = {Funkhouser, H. G.},
  journal   = {The American mathematical monthly},
  title     = {A short account of the history of symmetric functions of roots of equations},
  year      = {1930},
  number    = {7},
  pages     = {357--365},
  volume    = {37},
  publisher = {Taylor \& Francis},
}

@book{vinberg2003course,
  title={{A Course in Algebra}},
  author={Vinberg, {\.E}. B.},
  number={56},
  year={2003},
  publisher={American Mathematical Soc.}
}

@book{strang2022introduction,
  title={{Introduction to Linear Algebra}},
  edition = {Sixth Edition},
  author={Strang, G.},
  year={2022},
  publisher={SIAM}
}

@book{horn2012matrix,
  title={Matrix analysis},
  author={Horn, R. A. and Johnson, C. R.},
  year={2012},
  publisher={Cambridge university press}
}

@article{trefethen2011six,
  title={Six myths of polynomial interpolation and quadrature},
  author={Trefethen, L. N.},
  year={2011},
  publisher={Mathematics Today}
}

@Article{aurentz2015fast,
  author    = {Aurentz, J. L. and Mach, T. and Vandebril, R. and Watkins, D. S.},
  journal   = {SIAM J. Matrix Anal. Appl.},
  title     = {Fast and backward stable computation of roots of polynomials},
  year      = {2015},
  number    = {3},
  pages     = {942--973},
  volume    = {36},
  fjournal  = {SIAM Journal on Matrix Analysis and Applications},
  publisher = {SIAM},
}

@Article{aurentz2018fast,
  author    = {Aurentz, J. L. and Mach, T. and Robol, L. and Vandebril, R. and Watkins, D. S.},
  journal   = {SIAM J. Matrix Anal. Appl.},
  title     = {Fast and backward stable computation of roots of polynomials, Part II: Backward error analysis; companion matrix and companion pencil},
  year      = {2018},
  number    = {3},
  pages     = {1245--1269},
  volume    = {39},
  fjournal  = {SIAM Journal on Matrix Analysis and Applications},
  publisher = {SIAM},
}

@book{mac_duffee_theory_1933,
	address = {Berlin, Heidelberg},
	title = {The {Theory} of {Matrices}},
	isbn = {9783642984211},
	publisher = {Springer Berlin Heidelberg},
	author = {MacDuffee, C. C.},
	year = {1933},
	doi = {10.1007/978-3-642-99234-6},
}

@Article{loewy_begleitmatrizen_1920,
  author   = {Loewy, A.},
  journal  = {Math. Z.},
  title    = {Begleitmatrizen und lineare homogene {Differentialausdrücke}},
  year     = {1920},
  issn     = {1432-1823},
  month    = mar,
  number   = {1},
  pages    = {58--125},
  volume   = {7},
  doi      = {10.1007/BF01199397},
  fjournal = {Mathematische Zeitschrift},
  language = {de},
}

@Article{edelman1995polynomial,
  author   = {Edelman, A. and Murakami, H.},
  journal  = {Math. Comput.},
  title    = {Polynomial roots from companion matrix eigenvalues},
  year     = {1995},
  number   = {210},
  pages    = {763--776},
  volume   = {64},
  fjournal = {Mathematics of Computation},
}

@Article{mackey2013continuing,
  author    = {Mackey, D. S.},
  journal   = {Linear Algebra Appl.},
  title     = {The continuing influence of Fiedler’s work on companion matrices},
  year      = {2013},
  number    = {4},
  pages     = {810--817},
  volume    = {439},
  fjournal  = {Linear Algebra and its Applications},
  publisher = {Elsevier},
}

@article{Frobenius1879,
author = {Frobenius, G.},
journal = {Journal für die reine und angewandte Mathematik},
pages = {146-208},
title = {Theorie der linearen Formen mit ganzen Coefficienten.},
volume = {86},
year = {1879},
}

@Article{good1961colleague,
  author    = {Good, I.J.},
  journal   = {Q. J. Math.},
  title     = {The colleague matrix, a Chebyshev analogue of the companion matrix},
  year      = {1961},
  number    = {1},
  pages     = {61--68},
  volume    = {12},
  fjournal  = {The Quarterly Journal of Mathematics},
  publisher = {Oxford University Press},
}

@misc{DLMF,
         key = "{\relax DLMF}",
       title = "{\it NIST Digital Library of Mathematical Functions}",
edition = {Release 1.1.12 of 2023-12-15},
         url = "https://dlmf.nist.gov/",
        note = "F.~W.~J. Olver, A.~B. {Olde Daalhuis}, D.~W. Lozier, B.~I. Schneider,
                R.~F. Boisvert, C.~W. Clark, B.~R. Miller, B.~V. Saunders,
                H.~S. Cohl, and M.~A. McClain, eds."}

@book{rivlin_chebyshev_2020,
	address = {Mineola, New York},
	edition = {Second edition},
	title = {Chebyshev polynomials: from approximation theory to algebra \& number theory},
	isbn = {9780486842332},
	shorttitle = {Chebyshev polynomials},
	publisher = {Dover Publications, Inc},
	author = {Rivlin, T. J.},
	year = {2020},
}

@article{gautschi1972condition,
  title={The condition of orthogonal polynomials},
  author={Gautschi, W.},
  journal={Math. Comput.},
  volume={26},
  number={120},
  pages={923--924},
  year={1972}
}

@Article{gautschi1979condition,
  author   = {Gautschi, W.},
  journal  = {Math. Comput.},
  title    = {The condition of polynomials in power form},
  year     = {1979},
  number   = {145},
  pages    = {343--352},
  volume   = {33},
  fjournal = {Mathematics of Computation},
}

@Article{wilkinson1959evaluation1,
  author    = {Wilkinson, J. H.},
  journal   = {Numer. Math.},
  title     = {The evaluation of the zeros of ill-conditioned polynomials. Part I},
  year      = {1959},
  pages     = {150--166},
  volume    = {1},
  fjournal  = {Numerische Mathematik},
  publisher = {Springer},
}

@Article{wilkinson1959evaluation2,
  author    = {Wilkinson, J. H.},
  journal   = {Numer. Math.},
  title     = {The evaluation of the zeros of ill-conditioned polynomials. Part II},
  year      = {1959},
  pages     = {167--180},
  volume    = {1},
  fjournal  = {Numerische Mathematik},
  publisher = {Springer},
}

@Article{mosier1986root,
  author   = {Mosier, R. G.},
  journal  = {Math. Comput.},
  title    = {Root neighborhoods of a polynomial},
  year     = {1986},
  number   = {175},
  pages    = {265--273},
  volume   = {47},
  fjournal = {Mathematics of Computation},
}

@book{wilkinson2023rounding,
  title={Rounding errors in algebraic processes},
  author={Wilkinson, J. H.},
  year={2023},
  publisher={SIAM}
}

@book{trefethen2022numerical,
  title={Numerical linear algebra},
  author={Trefethen, L. N. and Bau, D.},
  volume={181},
  year={2022},
  publisher={Siam}
}

@Article{specht1956lage,
  author    = {Specht, W.},
  journal   = {Math. Nachr.},
  title     = {{Die Lage der Nullstellen eines Polynoms}},
  year      = {1956},
  number    = {5-6},
  pages     = {353--374},
  volume    = {15},
  fjournal  = {Mathematische Nachrichten},
  publisher = {Wiley Online Library},
}

@misc{serkh_provably_2021,
	title = {A {Provably} {Componentwise} {Backward} {Stable} $O(n^2)$ {QR} {Algorithm} for the {Diagonalization} of {Colleague} {Matrices}},
	doi = {10.48550/arXiv.2102.12186},
	publisher = {arXiv},
	author = {Serkh, K. and Rokhlin, V.},
	month = feb,
	year = {2021},
	note = {arXiv:2102.12186}
}

@book{trefethen2019approximation,
  title={Approximation theory and approximation practice},
edition = {Extended Edition},
  author={Trefethen, L. N.},
  year={2019},
  publisher={SIAM}
}

@Article{boyd2002computing,
  author    = {Boyd, J. P.},
  journal   = {SIAM J. Numer. Anal.},
  title     = {Computing zeros on a real interval through Chebyshev expansion and polynomial rootfinding},
  year      = {2002},
  number    = {5},
  pages     = {1666--1682},
  volume    = {40},
  fjournal  = {SIAM Journal on Numerical Analysis},
  publisher = {SIAM},
}

@Article{nakatsukasa2016stability,
  author   = {Nakatsukasa, Y. and Noferini, V.},
  journal  = {Math. Comput.},
  title    = {On the stability of computing polynomial roots via confederate linearizations},
  year     = {2016},
  number   = {301},
  pages    = {2391--2425},
  volume   = {85},
  fjournal = {Mathematics of Computation},
}

@Article{noferini2017chebyshev,
  author   = {Noferini, V. and P{\'e}rez, J.},
  journal  = {Math. Comput.},
  title    = {Chebyshev rootfinding via computing eigenvalues of colleague matrices: when is it stable?},
  year     = {2017},
  number   = {306},
  pages    = {1741--1767},
  volume   = {86},
  fjournal = {Mathematics of Computation},
}

@Article{cody1970survey,
  author    = {Cody, W.J.},
  journal   = {Siam Rev.},
  title     = {A survey of practical rational and polynomial approximation of functions},
  year      = {1970},
  number    = {3},
  pages     = {400--423},
  volume    = {12},
  fjournal  = {Siam Review},
  publisher = {SIAM},
}

@book{lanczos_applied_1988,
	address = {New York},
	title = {Applied analysis},
	isbn = {9780486656564},
	publisher = {Dover Publications},
	author = {Lanczos, C.},
	year = {1988},
}

@Article{thacher1964conversion,
  author    = {Thacher Jr, H. C.},
  journal   = {Commun. ACM},
  title     = {Conversion of a power to a series of Chebyshev polynomials},
  year      = {1964},
  number    = {3},
  pages     = {181--182},
  volume    = {7},
  fjournal  = {Communications of the ACM},
  publisher = {ACM New York, NY, USA},
}

@article{fox1968chebyshev,
  title={{Chebyshev polynomials in Numerical Analysis}},
  author={Fox, L. and Parker, I. B.},
  year={1968},
    location = {London},
  publisher={Oxford University Press}
}

@article{tucker2011validated,
  title={Validated numerics: a short introduction to rigorous computations},
  author={Tucker, W.},
  year={2011},
  publisher={Princeton University Press}
}

@book{mayer2017interval,
  title={Interval analysis: and automatic result verification},
  author={Mayer, G.},
  volume={65},
  year={2017},
  publisher={Walter de Gruyter GmbH \& Co KG}
}

@book{moore2009introduction,
  title={Introduction to interval analysis},
  author={Moore, R. E. and Kearfott, R. B. and Cloud, M. J.},
  year={2009},
  publisher={SIAM}
}

@Article{sitton2003factoring,
  author    = {Sitton, G. A. and Burrus, C. S. and Fox, J. W. and Treitel, S.},
  journal   = {IEEE Signal Process Mag.},
  title     = {Factoring very-high-degree polynomials},
  year      = {2003},
  number    = {6},
  pages     = {27--42},
  volume    = {20},
  fjournal  = {IEEE Signal Processing Magazine},
  publisher = {IEEE},
}

@article{noferini_structured_2021,
	title = {Structured backward errors in linearizations},
	volume = {54},
	issn = {1068-9613, 1068-9613},
	doi = {10.1553/etna_vol54s420},
	language = {en},
	journal = {ETNA - Electronic Transactions on Numerical Analysis},
	author = {Noferini, V. and Robol, L. and Vandebril, R.},
	year = {2021},
	pages = {420--442},
}

@article{neto2009electronic,
  title={The electronic properties of graphene},
  author={Castro Neto, A.H. and Guinea, F. and Peres, N.M.R. and Novoselov, K. S. and Geim, A. K.},
  journal={Rev. Mod. Phys.},
  volume={81},
  number={1},
  pages={109},
  year={2009},
  publisher={APS}
}

@article{wallace1947band,
  title={The band theory of graphite},
  author={Wallace, P. R.},
  journal={Physical review},
  volume={71},
  number={9},
  pages={622},
  year={1947},
  publisher={APS}
}

@Article{semenoff1984condensed,
  author    = {Semenoff, G. W.},
  journal   = {Phys. Rev. Lett.},
  title     = {Condensed-matter simulation of a three-dimensional anomaly},
  year      = {1984},
  number    = {26},
  pages     = {2449},
  volume    = {53},
  fjournal  = {Physical Review Letters},
  publisher = {APS},
}

@book{liu2018graphene,
  title={Graphene photonics},
  author={Liu, J.M. and Lin, I.T.},
  year={2018},
  publisher={Cambridge University Press}
}

@article{wang2023machine,
  title={Machine learning seams of conical intersection: A characteristic polynomial approach},
  author={Wang, T. Y. and Neville, S. P. and Schuurman, M. S.},
  journal={The Journal of Physical Chemistry Letters},
  volume={14},
  number={35},
  pages={7780--7786},
  year={2023},
  publisher={ACS Publications}
}

@Article{opalka2013interpolation,
  author    = {Opalka, D. and Domcke, W.},
  journal   = {J. Chem. Phys.},
  title     = {Interpolation of multi-sheeted multi-dimensional potential-energy surfaces via a linear optimization procedure},
  year      = {2013},
  number    = {22},
  volume    = {138},
  fjournal  = {The Journal of Chemical Physics},
  publisher = {AIP Publishing},
}

@Article{westermayr2020machine,
  author    = {Westermayr, J. and Marquetand, P.},
  journal   = {Chem. Rev.},
  title     = {Machine learning for electronically excited states of molecules},
  year      = {2020},
  number    = {16},
  pages     = {9873--9926},
  volume    = {121},
  fjournal  = {Chemical Reviews},
  publisher = {ACS Publications},
}

@Article{wilkinson2014excited,
  author    = {Wilkinson, I. and Boguslavskiy, A. E. and Mikosch, J. and Bertrand, J. B. and W{\"o}rner, H. J. and Villeneuve, D. M. and Spanner, M. and Patchkovskii, S. and Stolow, A.},
  journal   = {J. Chem. Phys.},
  title     = {Excited state dynamics in SO2. I. Bound state relaxation studied by time-resolved photoelectron-photoion coincidence spectroscopy},
  year      = {2014},
  number    = {20},
  volume    = {140},
  fjournal  = {The Journal of Chemical Physics},
  publisher = {AIP Publishing},
}

@Article{DrautzACE,
  title = {Atomic cluster expansion for accurate and transferable interatomic potentials},
  author = {Drautz, Ralf},
  journal = {Phys. Rev. B},
  volume = {99},
  issue = {1},
  pages = {014104},
  numpages = {15},
  year = {2019},
  publisher = {American Physical Society},
  doi = {10.1103/PhysRevB.99.014104},
}

@Article{ACECompleteness,
  author    = {Dusson, G. and Bachmayr, M. and Cs{\'a}nyi, G. and Drautz, R. and Etter, S. and van der Oord, C. and Ortner, C.},
  journal   = {J. Comput. Phys.},
  title     = {Atomic cluster expansion: Completeness, efficiency and stability},
  year      = {2022},
  pages     = {110946},
  volume    = {454},
  doi       = {10.1016/j.jcp.2022.110946},
  fjournal  = {Journal of Computational Physics},
  publisher = {Elsevier},
}

@InProceedings{Batatia2022mace,
  author    = {Batatia, I. and Kovacs, D. P. and Simm, G. and Ortner, C. and Cs{\'a}nyi, G.},
  booktitle = {Adv. Neural Inf. Process. Syst.},
  title     = {MACE: Higher Order Equivariant Message Passing Neural Networks for Fast and Accurate Force Fields},
  year      = {2022},
  editor    = {S. Koyejo and S. Mohamed and A. Agarwal and D. Belgrave and K. Cho and A. Oh},
  pages     = {11423--11436},
  publisher = {Curran Associates, Inc.},
  volume    = {35},
}

@article{kouppel1984multimode,
  title={Multimode molecular dynamics beyond the Born-Oppenheimer approximation},
  author={K{\"o}uppel, H. and Domcke, W. and Cederbaum, L. S.},
  journal={Advances in chemical physics},
  pages={59--246},
  year={1984},
  publisher={Wiley Online Library}
}

@Article{plasser2019highly,
  author    = {Plasser, F. and G{\'o}mez, S. and Menger, M. and Mai, S. and Gonz{\'a}lez, L.},
  journal   = {Phys. Chem. Chem. Phys.},
  title     = {Highly efficient surface hopping dynamics using a linear vibronic coupling model},
  year      = {2019},
  number    = {1},
  pages     = {57--69},
  volume    = {21},
  fjournal  = {Physical Chemistry Chemical Physics},
  publisher = {Royal Society of Chemistry},
}

@Article{axelrod_excited_2022,
  author    = {Axelrod, S. and Shakhnovich, E. and Gómez-Bombarelli, R.},
  journal   = {Nat. Commun.},
  title     = {Excited state non-adiabatic dynamics of large photoswitchable molecules using a chemically transferable machine learning potential},
  year      = {2022},
  issn      = {2041-1723},
  month     = jun,
  number    = {1},
  pages     = {3440},
  volume    = {13},
  copyright = {2022 The Author(s)},
  doi       = {10.1038/s41467-022-30999-w},
  fjournal  = {Nature Communications},
  language  = {en},
}

@Article{dral_molecular_2021,
  author    = {Dral, P. O. and Barbatti, M.},
  journal   = {Nat. Rev. Chem.},
  title     = {Molecular excited states through a machine learning lens},
  year      = {2021},
  issn      = {2397-3358},
  month     = jun,
  number    = {6},
  pages     = {388--405},
  volume    = {5},
  copyright = {2021 Springer Nature Limited},
  doi       = {10.1038/s41570-021-00278-1},
  fjournal  = {Nature Reviews Chemistry},
  language  = {en},
}

@Article{kulik_roadmap_2022,
  author   = {Kulik, H. J. and Hammerschmidt, T. and Schmidt, J. and Botti, S. and Marques, M. A. L. and Boley, M. and Scheffler, M. and Todorović, M. and Rinke, P. and Oses, C. and Smolyanyuk, A. and Curtarolo, S. and Tkatchenko, A. and Bartók, A. P. and Manzhos, S. and Ihara, M. and Carrington, T. and Behler, J. and Isayev, O. and Veit, M. and Grisafi, A. and Nigam, J. and Ceriotti, M. and Schütt, K. T. and Westermayr, J. and Gastegger, M. and Maurer, R. J. and Kalita, B. and Burke, K. and Nagai, R. and Akashi, R. and Sugino, O. and Hermann, J. and Noé, F. and Pilati, S. and Draxl, C. and Kuban, M. and Rigamonti, S. and Scheidgen, M. and Esters, M. and Hicks, D. and Toher, C. and Balachandran, P. V. and Tamblyn, I. and Whitelam, S. and Bellinger, C. and Ghiringhelli, L. M.},
  journal  = {Electron. Struct.},
  title    = {Roadmap on {Machine} learning in electronic structure},
  year     = {2022},
  issn     = {2516-1075},
  month    = aug,
  number   = {2},
  pages    = {023004},
  volume   = {4},
  doi      = {10.1088/2516-1075/ac572f},
  fjournal = {Electronic Structure},
  language = {en},
}

@Article{fedik_extending_2022,
  author    = {Fedik, N. and Zubatyuk, R. and Kulichenko, M. and Lubbers, N. and Smith, J. S. and Nebgen, B. and Messerly, R. and Li, Y.W. and Boldyrev, A. I. and Barros, K. and Isayev, O. and Tretiak, S.},
  journal   = {Nat. Rev. Chem.},
  title     = {Extending machine learning beyond interatomic potentials for predicting molecular properties},
  year      = {2022},
  issn      = {2397-3358},
  month     = sep,
  number    = {9},
  pages     = {653--672},
  volume    = {6},
  copyright = {2022 Springer Nature Limited},
  doi       = {10.1038/s41570-022-00416-3},
  fjournal  = {Nature Reviews Chemistry},
  language  = {en},
}

@Article{westermayr_deep_2022,
  author    = {Westermayr, J. and Gastegger, M. and Vörös, D. and Panzenboeck, L. and Joerg, F. and González, L. and Marquetand, P.},
  journal   = {Nat. Chem.},
  title     = {Deep learning study of tyrosine reveals that roaming can lead to photodamage},
  year      = {2022},
  issn      = {1755-4349},
  month     = aug,
  number    = {8},
  pages     = {914--919},
  volume    = {14},
  copyright = {2022 The Author(s), under exclusive licence to Springer Nature Limited},
  doi       = {10.1038/s41557-022-00950-z},
  fjournal  = {Nature Chemistry},
  language  = {en},
}

@article{xian_machine_2023,
	title = {A machine learning route between band mapping and band structure},
	volume = {3},
	issn = {2662-8457},
	doi = {10.1038/s43588-022-00382-2},
	language = {en},
	number = {1},
	journal = {Nature Computational Science},
	author = {Xian, R. P. and Stimper, V. and Zacharias, M. and Dendzik, M. and Dong, S. and Beaulieu, S. and Schölkopf, B. and Wolf, M. and Rettig, L. and Carbogno, C. and Bauer, S. and Ernstorfer, R.},
	month = jan,
	year = {2023},
	pages = {101--114},
}

@book{phani_dynamics_2017,
	address = {Chichester, West Sussex, United Kingdom},
	title = {Dynamics of lattice materials},
	isbn = {9781118729595},
	publisher = {John Wiley \& Sons, Inc},
	editor = {Phani, A. S. and Hussein, M. I.},
	year = {2017},
}

@Article{sadat_machine_2020,
  author   = {Sadat, Seid M. and Wang, Robert Y.},
  journal  = {J. Appl. Phys.},
  title    = {A machine learning based approach for phononic crystal property discovery},
  year     = {2020},
  issn     = {0021-8979, 1089-7550},
  month    = jul,
  number   = {2},
  pages    = {025106},
  volume   = {128},
  doi      = {10.1063/5.0006153},
  fjournal = {Journal of Applied Physics},
  language = {en},
}

@Article{westermayr_combining_2020,
  author     = {Westermayr, J. and Gastegger, M. and Marquetand, P.},
  journal    = {J. Phys. Chem. Lett.},
  title      = {Combining {SchNet} and {SHARC}: {The} {SchNarc} {Machine} {Learning} Approach for {Excited}-{State} {Dynamics}},
  year       = {2020},
  issn       = {1948-7185},
  month      = may,
  number     = {10},
  pages      = {3828--3834},
  volume     = {11},
  doi        = {10.1021/acs.jpclett.0c00527},
  fjournal   = {Journal of Physical Chemistry Letters},
  pmcid      = {PMC7246974},
  pmid       = {32311258},
  shorttitle = {Combining {SchNet} and {SHARC}},
}

@Article{mai_molecular_2020,
  author     = {Mai, Sebastian and González, Leticia},
  journal    = {Angew. Chem. Int. Ed.},
  title      = {Molecular {Photochemistry}: {Recent} {Developments} in {Theory}},
  year       = {2020},
  issn       = {1433-7851, 1521-3773},
  month      = sep,
  number     = {39},
  pages      = {16832--16846},
  volume     = {59},
  doi        = {10.1002/anie.201916381},
  fjournal   = {Angewandte Chemie International Edition},
  language   = {en},
  shorttitle = {Molecular {Photochemistry}},
}

\end{document}